\documentclass[a4paper,11pt,english]{article} 
\usepackage{amsmath}
\usepackage{amsfonts}
\usepackage{amssymb}
\usepackage[all]{xypic}
\usepackage{graphics}

\usepackage{babel}
\usepackage{hyperref}
\usepackage{here}
\usepackage{verbatim}

\def\a{\alpha}
\def\b{\beta}

\def\vp{\varphi}
\def\g{\gamma}

\def\o{\omega}

\def\nbh{neighborhood}

\def\bP{{\mathbb P}}

\def\bR{{\mathbb R}}
\def\bZ{{\mathbb Z}}
\def\bC{{\mathbb C}}

\def\b1{{\rm id}}

\newfont{\goth}{eufm10 scaled \magstep1}
\def\ga{\mbox{\goth a}}

\def\gc{\mbox{\goth c}}

\def\gg{\mbox{\goth g}}
\def\gh{{\mbox{\goth h}}}
\def\gk{\mbox{\goth k}}
\def\gl{\mbox{\goth l}}
\def\gm{\mbox{\goth m}}
\def\gn{\mbox{\goth n}}

\def\gt{\mbox{\goth t}}

\newfont{\mcal}{eusm10 scaled \magstep1}

\def\ch{\mbox{\mcal H}}

\def\cl{\mbox{\mcal L}}

\newtheorem{Th}{Theorem}
\newtheorem{theorem}{Theorem}
\newtheorem{proposition}[theorem]{Proposition}
\newtheorem{definition}[theorem]{Definition}
\newtheorem{lemma}[theorem]{Lemma}
\newtheorem{corollary}[theorem]{Corollary}
\newtheorem{remark}[theorem]{Remark}
\newtheorem{remar}[theorem]{Remark}
\newtheorem{Prop}[Th]{Proposition}
\newtheorem{Cor}[Th]{Corollary}
\newtheorem{Lem}[Th]{Lemma}
\newtheorem{Def}[Th]{Definition}
\newtheorem{Ex}[Th]{Example }
\def\bt{\begin{Th}}
\def\et{\end{Th}}
\def\bp{\begin{Prop}}
\def\ep{\end{Prop}}
\def\bc{\begin{Cor}}
\def\ec{\end{Cor}}
\def\bl{\begin{Lem}}
\def\el{\end{Lem}}
\def\bd{\begin{Def}}
\def\ed{\end{Def}}
\def\bex{\begin{Ex}}
\def\eex{\end{Ex}}
\def\br{\begin{remar}}
\def\er{\end{remar}}
\def\pf{\noindent{\it Proof:\ }}
\def\qed{\hfill$\square$}

\def\be{\begin{equation}}
\def\ee{\end{equation}}
\def\ben{\begin{enumerate}}
 \def\een{\end{enumerate}}
\def\ba{\begin{array}{rlll}}
\def\ea{\end{array}}
\def\bea{\begin{eqnarray}}
\def\eea{\end{eqnarray}}
\def\bean{\begin{eqnarray*}}
\def\eean{\end{eqnarray*}}

\def\ad{{\rm ad}\,}
\def\Ad{{\rm Ad}\,}

\def\tr{\mathrm{tr\;}}

\def\ker{\mathrm{ker\;}}
\def\diag{\mathrm{diag\;}}

\def\C1{cohomogeneity one  }
\def\Km{K\"ahler manifold \,}

\begin{document}

\title{Cohomogeneity one K\"ahler  and K\"ahler-Einstein manifolds  with   one  singular orbit I }
\author{Dmitri Alekseevsky\footnote{A.A.Kharkevich Institute for Information Transition
Problems, B.Karetnuj  per.,19, 127051, Moscow, Russia and
Faculty of Science  University of Hradec Kralove,  Rokitanskeho 62, Hradec Kralove,
50003, Czech Republic } \and  Fabio Zuddas\footnote{Dipartimento di Scienze Matematiche, Informatiche e Fisiche, Via delle Scienze 206, Udine (Italy)}}



\maketitle

 \abstract
Let  $M$  be a \C1  manifold of a  compact semisimple Lie group $G$  with one  singular orbit $S_0 = G/H$.  Then $M$ is $G$-diffeomorphic
 to the   total space  $G \times_H V$ of  the  homogeneous vector  bundle  over $S_0$ defined  by a sphere transitive  representation of
 $G$ in a vector  space $V$.  We describe  all  such manifolds $M$ which admit  an invariant K\"ahler  structure  of  standard  type. This means  that
 the  restriction  $\mu : S = Gx =  G/L \to  F = G/K $ of the  moment map of $M$  to a regular orbit $S=G/L$ is a  holomorphic  map
of  $S$  with the induced  CR  structure  onto  a   flag manifold $F = G/K$, where  $K = N_G(L)$, endowed with an invariant   complex  structure  $J^F$.
We  describe all such  standard K\"ahler  \C1 manifolds  in  terms of the painted  Dynkin diagram  associated  with $(F = G/K,J^F)$   and a parametrized interval in  some  $T$-Weyl chamber.   \\
We determine   which of these  manifolds  admit invariant K\"ahler-Einstein metrics.

\tableofcontents

\section{Introduction and statement of the results}

We  will study \C1 K\"ahler  $G$-manifolds  of a  compact  semisimple Lie  group $G$.
By a \C1  manifold  we understand  an $n$-dimensional manifold  $M$  together  with  a proper  action of a connected  Lie group $G$ which has a (real) codimension one  orbit.\\
It  is called  a Riemannian (respectively, complex; K\"ahler ) \C1  manifold  if  an invariant Riemannian  metric
 $g $ (respectively, an  invariant  complex  structure $J$;  an invariant  K\"ahler   structure   $(J, \o)$) is given,  where  $J$ is  a  complex  structure   and $\o$    is a  symplectic     form    such   that  $g := - \o \circ J = \o(\cdot, J \cdot)$   is  a  K\"ahler metric.\\

Following  \cite{P-S} we   will    consider K\"ahler  \C1 $G$-manifolds $(M, J, \o)$  of  the { \it standard  type},
  that is  manifolds   which satisfy  the  following   conditions:\\
(i)   Any regular  orbit $ S= G x = G/L$ is  an {\it ordinary} manifold. This
means that the  normalizer   $K = N_G(L)$ of  the  stability
subgroup is  the centralizer  of a  torus in $G$ and $\dim K/L =1$.\\
 (ii) the CR   structure $(\ch, J^{\ch})$  induced  by the complex
structure $J$ of $M$  on  a (codimension one) regular orbit  $S = G/L$  is {\it projectable},
that is the projection  $\pi: S=G/L \to F=G/K$ is a holomorphic map
of a CR manifold  onto the flag manifold  $F$ equipped  with a
fixed invariant complex  structure $J^F$   which  does not  depend on $S$.\\
Condition (ii) depends on the complex  structure  $J$ on $M$    and shows   that the  CR   structure on  a  regular orbit  $G/L$ is  determined  by a   fixed invariant  complex  structure $J_F$ on   the  flag manifold $F$.  In particular, all  regular orbits  are isomorphic  as   homogeneous   CR manifolds.\\
  The  conditions (i),(ii)  imply  that    the  moment map $\mu : M \to \gg^* \simeq \gg$ of  the  symplectic
 $G$-manifold $(M, \o)$ maps  any regular orbit  $(S= Gx = G/L, J^{\ch})$ holomorphically to  the  same  flag manifold  $(F = G/K, J^F)$.
  Note  that   $\pi:= \mu|_S : S = G/L \to  F = G/K $ is the  natural equivariant projection.\\
\smallskip

A  homogeneous CR  manifold  $(S = G/L, \mathcal{H}, J^{\mathcal{H}})$ which satisfies conditions
(i),(ii)  is called  a  {\it  standard  homogeneous CR manifold}.\\
 So, we can equivalently say that a  complex (in particular,  K\"ahler) \C1 manifold is  of  the  {\it  standard  type} { or, shortly, {\it  standard }
if all  regular orbits  are   standard CR manifolds  associated  with  a  fixed flag manifold $(F=G/K,J^F)$  with   a  complex  structure $J^F$.\\

The condition  (i) is  a weak condition. It
 is  equivalent  to  the  condition that   the Lie  algebra $\gl = \mathrm{Lie}\, L$ is not   the centralizer of  a regular 3-dimensional subalgebra  of the Lie  algebra  $\gg = \mathrm{Lie}\, G $.\\
 The classification of  all non standard  homogeneous CR manifolds  with non degenerate Levi form  is  given in \cite{A-S}.\\
\smallskip

Investigation of invariant  Einstein metrics on cohomogeneity one manifolds had been   started   by  D. Page  and C.Pope in \cite{P-P},   where  they  construct the first example of  such metrics.   L.Berard- Bergery  in  his  famous paper   (\cite{Ber})  developed
a  systematic approach  for  construction of  invariant Einstein metrics on \C1  manifolds.   A deep investigation  of  singular ODE  for invariant   Einstein metrics on \C1 manifolds  had been   done  by  Eschenburg and Wang \cite{E-W} .\\

 In the K\"ahler case, Y. Sakane  gives in  \cite{S} conditions for the existence of K\"ahler-Einstein metrics on $\bC\bP^1$-bundles  $P$  over Hermitian symmetric spaces of compact type which are of cohomogenity one with respect to a maximal compact subgroup of the automorphism group of $P$. In \cite{K-S1} and \cite{K-S2} more examples are obtained from these $\bC\bP^1$-bundles by blowing down. We refer the reader also to \cite{H-S}, (see also \cite{B} and the references therein).\\

Invariant  K\"ahler-Einstein metrics  on    cohomogeneity one manifolds $M$ of a  compact  semisimple Lie group $G$
 had been  studied in  two important papers  by  F. Podest\'a  and  A. Spiro \cite{P-S}  and  A. Dancer   and  M. Y.  Wang  \cite{D-W}.  In both papers    the  authors  reduced   the  K\"ahler-Einstein
 equation   for  an invariant  metric  in  the   regular  open  submanifold   $M_{reg} = G/L \times (a,b)$  of  $M$ to an ODE   for  one  function   together   with    some algebraic    conditions,  which  are  described in terms  of  the  reductive decomposition  associated  with a  regular orbit  $G/L$.
 They  solved    this   equation   and   find necessary and  sufficient  conditions for   extendibility  of the K\"ahler-Einstein metric in  $M_{reg}$   to the whole manifold $M$. They considered some examples of manifolds which satisfy these conditions, but did not study such manifolds systematically.\\

The main aim of this paper is  to  give  a description of  (non
compact) standard \C1 K\"ahler and K\"ahler-Einstein  manifolds
with one singular orbit in terms of  painted  Dynkin diagrams.
  We   closely follow the approach   by F. Podest\`a and A. Spiro
\cite{P-S},  who give  a useful description of  standard  \C1
K\"ahler   manifolds $M$  in terms of  ``abstract models"  and  get
an effective criterion of existence of an invariant  K\"ahler-Einstein
metric  on a (compact)  \C1   manifold with two singular  orbits.
We   will reformulate   the basic results of  \cite{P-S}
  and simplify the proofs.\\

The  structure of the paper is  the  following. In the Preliminaries  we recall the  basic  facts  about \C1 Riemannian manifolds, CR manifolds and flag manifolds which  we  need. In Section 3  we    define     standard  homogeneous  CR  manifolds and  standard \C1  complex  and K\"ahler manifolds  and  discuss  their  properties.\\
   A standard    homogeneous  CR  manifold   $S = G/L$  with  associated  flag manifold  $F = G/K = G/N_G(L)$ is defined  by  the    standard  (reductive)    decomposition
$$   \gg   =  \gl   + \bR Z_F^0 + \gm = \gk + \gm$$
orthogonal  with respect to the  Killing  form $B$  of  $\mathfrak{g}$, where $\gl = \mathrm{Lie}\, L,   \,  \gk = \mathrm{Lie}\, K
   = \gl + \bR Z_F^0$ and $Z_F^0$  is an $\Ad_L$-invariant
 vector (called {\it fundamental vector}), normalized  by $B(Z_F^0,Z_F^0)=-1$.
 We    will identify   $Z_F^0$  with  a  $G$-invariant vector  field on  $S$ which is  the   fundamental vector   field   of the
 principal  $T^1$-bundle  $\pi :S=G/L \to  S/T^1 =  G/K = G/L \cdot T^1$.
  An $\Ad_K$-invariant  (integrable)  complex  structure
  $J^{\gm}$ in $\gm$
     defines    a  complex  structure  $J_F$
in the  flag manifold $F = G/K$  with the reductive  decomposition
$\gg = \gk + \gm$ and an invariant projectable CR structure  $ (\ch,
J^{\ch})$ on $S = G/L$
where $\ch$ is    the invariant  distribution defined  by the   subspace  $\gm$.\\

 In Section 4 we  describe  standard  invariant K\"ahler   structures  on a   regular  \C1  manifold
  $M_{reg} = (0,d) \times  G/L$ . Following \cite{P-S}, we  show  that  such    structures  are  in  1-1 correspondence
  with the    parametrized  open  intervals  $(Z_0Z_d) \subset C(J^F)$
  which are  parallel to $Z^0 = -i Z^0_F$.
   Here  $C(J^F)$  denotes the $T$-Weyl chamber (i.e. roughly speaking, the projection of  the Weyl  chamber of $\gg$ to the center
   $Z(\gk)$  of $\gk$)  associated   with  the  invariant complex   structure $J^F$
     of  the  flag  manifold $F = \mu(S)= G/K $. \\
 In Section 5 we  give a new proof  of Podest\`a-Spiro   formula  for  the Ricci form $\rho$.  The Einstein  equation
   for an invariant K\"ahler metric  in $M_{reg}$  reduces  to a second order ODE   for the
     function  $f(t)$  which   defines  a parametrization  of  the   interval  $(Z_0 Z_d) \subset C(J^F)$   and  a linear
     relation   between the  initial vector $Z_0$, the vector $Z^0$ and  the  Koszul vector
      $Z^{Kos} \in C(J^F)$  which defines  the  invariant K\"ahler-Einstein metric on the  flag manifold  $F =
      (G/K, J_F)$. \\
      To  calculate $\rho$,  we construct holomorphic  coordinates  $\{z_i\}$   in
$M_{reg} = (0,d) \times  G/L$ which are
 an extension of  local  holomorphic  coordinates in the  flag manifold $F = G/K = G/L \cdot T^1$ ,   and  use  the  formula $\rho = -i \partial \bar \partial \log \mu$,  where $\mu(z_i, \bar z_i)$ is  the  density  associated  with the  volume  form
 $$vol = \mu(z_i, \bar z_i) dz^m \wedge d \bar z^m,\,   m = \dim_{\bC}M_{reg}.$$
The    singular  orbit $S_0$ of a  standard  \C1   manifold is   a  complex  submanifold, hence  a   flag manifold $(S_0= G/H, J_S)$.  In Section 6,  we    describe   all   standard \C1  manifolds   with     fixed   singular orbit   $(S_0 = G/H,J_S)$. Any  such manifold
$M= M_{\varphi}$  is  defined   by  a   surjective   homomorphism $\varphi : H \to U_m $   and is   the  total   space  of   the homogeneous  vector  bundle  $M_{\varphi} = G \times_H \mathbb{C}^m_{\varphi} \to S_0 =G/H$      defined  by  $\varphi$.
  The  flag manifold  $(S_0,J_S)$    is  determined    by   a  painted Dynkin   diagram   (PDD).  In  terms  of PDD,   the   homomorphism  $\varphi$ is defined  by  a connected component (a  white string) of   type  $A_{m-1}$ of  the  white   subdiagram of  PDD  and   a  character   $\chi : Z(H) =T^k  \to  T^1$ of  the  center $Z(H)$. The  complex  structure  $J$ on $M_{\varphi}$ is  the natural  extension of  the  complex  structure
  $J_F$.  If $e \neq 0$ is  a  vector  from  $\mathbb{C}^m$, then   regular  orbits $S_t :=  G \times_H (te) = G/L$   are  parametrized   by
   $t>0$  where   $L = H_{te}$ is  the   stabilizer. The     subgroup $K =  H_{[e]}$  is  the  stabilizer  of  the  line $[e] \in P \mathbb{C}^m$.
   This  shows   that $\varphi$   determines  the  standard  reductive   decomposition
   $$    \mathfrak{g} = \mathfrak{l} + \mathbb{R}Z_F^0 + \mathfrak{m}.                 $$

An invariant  K\"ahler metric in $M_{\varphi}$ is  defined  by   an interval  $(Z_0 Z_d) \subset  C(J_F)$  of    the  $T$-Weyl  chamber associated   with   $J_F$  which   starts   from a   face associated    with $G/H$     together   with   a parametrization  which  satisfies   the Verdiani  boundary  condition. \\
In  the  last  chapter   we   give necessary   and  sufficient  conditions  for  a manifold $M_{\varphi}$  to have  an invariant K\"ahler-Einstein metric.  We  will use  this  condition in  the   second part of  this paper  for  explicit   description of   such K\"ahler-Einstein metrics on standard  \C1   manifolds  of  a classical Lie  group  $G$.\\

\medskip

\noindent {\bf Acknowledgments}. The authors would like to thank
F. Podest\`a and A. Spiro for useful discussions,  suggestions and clarifications on their   results.

\section{Preliminaries}
\subsection{Riemannian cohomogeneity one manifolds}
Let $G$ be a   compact Lie  group  and $(M,g)$ a Riemannian \C1  $G$-manifold, that  is  $G$ is  an isometry group of  $(M,g) $
 with a (real)  codimension one regular orbit  $S = Gx = G/L$. Denote  by  $\pi : M \to \Omega = M/G$ the natural projection   to the orbit  space.
 There   are  four   cases: the  orbit   space is  diffeomorphic to a) (0,1) ,  b) [0,1) c)
[0,1] or d) $S^1$.\\
 In the  case b)  there is   one  singular  orbit $\pi^{-1}(0) = G/H_0$   and   in the case   c)
   there  are two singular orbits   $S_{\epsilon} = \pi^{-1}(\epsilon) = G/H_{\epsilon},\,\, \epsilon = 0,1$.
A naturally parametrized geodesic $\gamma(t)$ normal  to an orbit
remains orthogonal to  any orbit  and it is called  a { \it normal
geodesic}.
 If it is  complete, it intersects   any orbit. In the cases b), c), we  will assume that $x= \gamma(0) $ belongs
to the  singular orbit $S_0 = G x = G/H_0$. Then the stabilizer
$H_0$ transforms  $\gamma$ to any other normal geodesic through
$x$ and the  isotropy  representation $j(H_0)|V$ restricted  to   the
normal space $V = T_x^{\perp}(S_0)$ acts transitively  on the
sphere: in other  words,   the  orbit $j(H)v = H/H_v = H/L$ is  the
sphere.   The \C1 $G$-manifold $M$ near $S_0$ is locally
$G$-diffeomorphic to the total space $G \times_H V $ of a
homogeneous  vector bundle over the singular orbit $S_0$.\\
In the case c), when     $M/G \simeq [0,1]$, rescaling the metric we
may assume  that $\gamma(\epsilon) \in S_{\epsilon},\,\, \epsilon =
0,1 $. Then   the \C1 manifold $M$ is determined   by the triple
$(H_0, L, H_1)$ of stability subgroups of  $\gamma(0), \gamma(1/2),
\gamma(1)$ and is denoted by $M(H_0,L,H_1)$. Note that $L \subset
H_0 \cap H_1$  and $H_{\epsilon}/L $ are  spheres. \\
In the  case d), $M/G \simeq S^1 = \{   exp (2 \pi i t) \}$,  we
assume that $\gamma(0), \gamma(1)
 \in  S_0 = \pi^{-1} (0)$.
Note  that the  stabilizer   $L = G_{\gamma(t)}$   of  any  regular point  of  $\gamma$   preserves
  pointwisely  the geodesic $\gamma$. So  we can identify
  any  regular  orbit   $S_t = G (\gamma(t) )$  with the  same  homogeneous  space  $G/L$.
    Deleting   singular  points  (if  they exist)   or,  in the
    case d),  the regular  orbit $S_0 = \pi^{-1}(0)$   we   get an open dense   submanifold  $M_{reg}$  of
      regular points  which is  $G$-diffeomorphic  to   $M_{reg} = (0,1) \times  G/L$.

Note   that  the orbit   space $\Omega$ of  a Riemannian \C1
manifold  has  a structure of a metric  space.
 The  following proposition  gives a description of \C1 Riemannian manifolds
 $(M,g)$ and  their orbit  spaces. 

\begin{proposition}\label{RiemC1manifolds} Let $(M,g)$  be a Riemannian $G$-manifold with the orbit  space $\Omega$. Then up to a  homothety, $(M,g)$ is described as follows:\

a)(No  singular  orbit)  \,  $M = \Omega \times G/L$,  where
$\Omega \approx (0,1)$. Moreover, if a normal geodesic  is complete,
then $\Omega = \mathbb{R}$. In the non-complete case, $\Omega =
(0,1) $ or $ \bR^+$. The  metric  is  given  by
$$  g =  dt^2 + g_t$$
 where
 $g_t , \, t \in \Omega$ is  a 1-parameter family of invariant Riemannian metrics on $G/L$.\\

 b)(One  singular orbit   $S_0 = G/H$)\, The  orbit  space  is $\Omega =  [0,d),\,\,  d= \infty$ or
 $1$. If  a normal geodesic  is  complete, then  $d=\infty$    and  the  manifold
$M = M(H,L)= G \times_H V$ is the  homogeneous vector  bundle  over
  the  singular orbit  $S_0$ defined  by a  sphere  transitive  orthogonal  representation
   $\nu : H \to O(V)$  of $H$  in an Euclidean vector  space $V$.
   In  the non complete case $M$ is
    a tubular  invariant \nbh \,  $M = G \times_H B \subset M(H,L)$   of the zero section
    where $B$  is the unit ball in $V$.
   The  invariant metric in $M$  is  an  invariant extension  of   the $j(H)$-invariant   Riemannian metric
    in $V$, which is  described  by L. Verdiani \cite{V}.\\

c)(Two singular orbits $S_{\epsilon} = G/H_{\epsilon},\,\, \epsilon
= 0,1$)\,, $\Omega = [0,1]$. The  Riemannian manifold $M$ is
obtained by gluing together two manifolds $M_- = \pi^{-1}([0,1/2)),
\, M_+ = \pi^{-1}( ( 1/2,1]) $ of  type  b)  along the isomorphic
boundary $\partial M_{\pm} = G/L $.    As a \C1 manifold it is
defined  by   the  triple of  subgroups   $H_0 ,L, H_1$  such that
$H_{\epsilon}/L = S^{n_{\epsilon}}$ and it is denoted by   $M = M(H_0,L,H_1)$.\\

d) (No singular orbit, $\Omega $ is not  contractible) \,\,
 $\Omega = S^1$ and  $M$ is a fibre bundle over   the  circle $S^1$  having as
  universal cover  a Riemannian manifold $\bR \times  G/L$ of  type a).\\
\end{proposition}

\subsection{ Flag manifolds and painted Dynkin diagrams}

\subsubsection{Isotropy decomposition, $T$-roots, $T$-Weyl chambers
  and  invariant complex  structures}

\noindent  Let $F = G/K = \mathrm{Ad}_G Z$, where $Z \in
\mathfrak{g}$, be a  flag manifold, i.e. an adjoint orbit of a
 compact semisimple Lie group $G$
  with   the  $B$-orthogonal \\
   (where $B$ is  the Killing form) reductive   decomposition
  $$ \mathfrak{g} = \gk + \gm  = C_{\gg}(Z) + \gm.$$
   We can decompose  $\gk$  as
     $$\gk   = Z(\gk) \oplus \gk'$$
   where  $\gk'$ is  the semisimple part   and  $Z(\gk)$ is  the center.
    We fix a Cartan subalgebra  $\gc$ of  $\gk$ (hence  also of
    $\gg$) and   denote  by
      $R$  the   root  system
     of   the   complex Lie  algebra $\gg^{\mathbb{C}}$ w.r.t. the Cartan  subalgebra
      $\gc^{\mathbb{C}}$.  We set
      $$R_{\mathfrak{k}} :=\{ \a \in R, \, \a(Z(\gk)) = 0 \},\,R_{\mathfrak{m}} := R \setminus R_{\mathfrak{k}}. $$

 Then
  $$ \gk = \gc  + \gg (R_{\mathfrak{k}})^{\tau},\,   \gm = \gg(R_{\mathfrak{m}})^{\tau}, $$
where for a  subset $P \subset R$, we set
$$\gg(P) = \sum_{\alpha \in
P} \gg_{\alpha}$$

\noindent being $\gg_{\a}$ the  root  space with root  $\a$ and
$V^{\tau}$ means the fix point set in $V\subset \gg^{\bC}$  of the
complex conjugation $\tau$. Recall that the Killing   form induces
an Euclidean metric in the real vector space $i \gc $ and roots are
identified  with real linear forms on $i \gc$. We  set $\gt:= i
Z(\gk) \subset i \gc$ and denote  by
 $$\kappa: R \to R|_{\gt},\,\,  \a \mapsto \bar \a := \a|_{\gt}$$
  the  restriction map.

\begin{definition} 
The set $R_T =\kappa(R_{\mathfrak{m}})= R_{\mathfrak{m}}|{\gt}$ of linear forms
on $\gt$ which are  restriction of  roots from $R_{\mathfrak{m}}$ is called the system of
{\bf $T$-roots} and connected components $C$ of the set $\gt
\setminus \{ \ker \bar \a ,\, \bar \a \in R_T \}$
  are called  {\bf $T$-Weyl chambers}.
\end{definition}

Sets of  $T$-roots $\xi$  bijectively correspond to irreducible   $\gk$-submodules
 $ \gm(\xi):=  \gg(\kappa^{-1}(\xi))$ of the   complexified  isotropy module
 $\gm^{\mathbb{C}}$ of  the  flag manifold $F =G/K$.

So a decomposition of the $\gk$-modules  $\gm^{\mathbb{C}}$ and $\gm$
into irreducible
  submodules   can be  written
 as
$$  \gm^{\mathbb{C}} = \sum_{\xi \in R_T} \gm(\xi),\,\,
  \gm = \sum_{\xi \in R^+_T}[ \gm(\xi) + \gm(-\xi) ]^{\tau}$$
 where $R^+_T :=\kappa(R_{\mathfrak{m}}^+)$ is   the  system of  positive $T$-roots associated    with  a
system of positive  roots $R^+$, see  \cite{A-P}, \cite{A}.

We  fix a system of  simple  roots $\Pi_W$ of  $R_{\mathfrak{k}}$
and denote  by  $\Pi = \Pi_W \cup \Pi_B$ its extension  to a system
of simple roots of $R$.  Let $R^+ = R^+(\Pi) $  be the associated
system of positive  roots and $R^+_{\mathfrak{m}}:= R^+ \cap
R_{\mathfrak{m}}$. The  set  $R^+_T:= \kappa(R^+_{\mathfrak{m}})$ is
called positive $T$-root  set.

 We need  the following

\begin{theorem}\cite{A-P} 
There exists a one-to-one correspondence  between
extensions $\Pi = \Pi_W \cup \Pi_B$  of  the system $\Pi_W$ of
simple system of  $R_{\mathfrak{k}}$,  $T$-Weyl chambers $C\subset
\gt$ and invariant complex structures (ICS) $J$  on $F = G/K$. If
$\Pi_B = \{ \beta_1, \dots ,\beta_k \}$, then  the corresponding
$T$-Weyl chamber  is defined  by  $   C = \{   \bar \beta_1>0,
\dots ,\bar \beta_k >0 \}$ where $\bar \beta = \kappa(\beta)$  and
the  complex  structure is  defined  by $\pm i$-eigenspace decomposition
 \be \label{InvCompStrin"m"} \gm^{\bC} = \gm^+ + \gm^- =
\gg(R^+_{\mathfrak{m}}) + \gg(-R^+_{\mathfrak{m}}) \ee
 of the
complexified  tangent  space  $\gm^{\bC} = T_{eK}(G/K)$.
 \end{theorem}

The extension $\Pi = \Pi_W \cup \Pi_B$  can be graphically described
by  a painted Dynkin diagram, i.e.  the  Dynkin  diagram  which
represents the  system $\Pi$  with the nodes representing  $\Pi_B$
painted in  black. Such a diagram, which we   sometimes  identify
with  the pair $(\Pi_W, \Pi_B)$,  allows  to reconstruct   the flag
manifold $F= G/K$ with invariant  complex  structure $J^F$ as
follows: the  semisimple  part $\gk'$ of the  (connected) stability
subalgebra $\gk$ is  defined as the  regular semisimple  subalgebra
associated   with the closed  subsystem $R_{\mathfrak{k}}= R \cap
\mathrm{span}(\Pi_W)$ and the vectors $ih_{j}$ defined  by condition
$$ \beta_k(h_j) = \delta_{kj}, \, \alpha_i(h_j)=0,\,  \beta_j \in \Pi_B, \alpha_i \in \Pi_W$$
 form a basis of
the center $Z(\mathfrak{k})$. The complex structure is defined  by
(\ref{InvCompStrin"m"}).

\subsubsection{Invariant symplectic  forms  and  K\"ahler  structures }
An element $Z \in \mathfrak{t}$  is called   to  be {\bf
$K$-regular} if  its   centralizer  $C_G(Z) = K$ or, equivalently,
any $T$-root has a non-zero value on $Z$.

\begin{proposition}\label{propcorrespmetrics} 
(\cite{B-H}, \cite{A-P}) There exists
a natural one-to-one correspondence between elements $Z \in
\mathfrak{t}$ and closed invariant 2-forms  $\omega_Z$ on $G/K$,
given by
$$Z \leftrightarrow \omega_Z|_o = i \, d(B \circ Z),$$
\noindent where $d$
is  the exterior differential in  the Lie  algebra $\mathfrak{g}$
defined  by $d\alpha(X,Y) = - 1/2\alpha([X,Y])$  and  $o = eK \in G/K$.\\
Moreover,  regular elements $Z \in C$  from a $T$-Weyl chamber  $C$
correspond to the  K\"ahler  forms  $\omega_Z$ with respect to the
complex  structure   $J(C)$ associated to $C$, that is they define an invariant  K\"ahler
structure  $(\omega_Z, J(C))$.  The  2-form $ \frac{1}{2\pi}\omega_Z$
is integral if  the  1-form  $B\circ Z$  has integer  coordinates
 with respect  to  the  fundamental weights  $ \pi_i$   associated  with
  the  system  of black simple roots $ \beta_i \in \Pi_B$. 
\end{proposition}

\noindent Recall that if $\Pi_W = \{ \a_1, \dots, \a_m\}$ (resp. $\Pi_B = \{ \beta_1, \dots, \beta_k\}$) is the set of
 white
(resp. black)  simple roots, then the fundamental weight $ \pi_i$
associated with $\beta_i$, $i = 1, \dots, k$, is the linear form
defined by

\begin{equation}\label{deffundwe}
\frac{2\langle  \pi_i, \beta_j \rangle}{\| \beta_j \|^2} =
\delta_{ij}, \ \ \langle  \pi_i, \alpha_j \rangle = 0.
\end{equation}

\noindent The $B$-dual  to $ \pi_i$ vectors $h_i$ form
  a basis of $\gt$.

\noindent Let $E_{\alpha} \in \gg_{\alpha}, \, \alpha \in R$, be the {\it Chevalley basis} of $\mathfrak{g}(R)$
such that $B(E_{\alpha}, E_{- \alpha}) = \frac{2}{<\alpha, \alpha>}$
where $< .,.>$ is the scalar product  in $i\gc^* = \mathrm{span}(R)$
induced by the Killing  form.  We denote by   $\omega_{\alpha} =
B\circ E_{\alpha}$ the dual basis of 1-forms. Then for $Z \in \gt$

\be \label{omega_Z}
 \omega_Z = -i\sum_{\alpha \in R_{\mathfrak{m}}^+}
\frac{2 \alpha(Z)}{<\alpha, \alpha>} \omega_{\alpha} \wedge
\omega_{-\alpha}\ee

\noindent Indeed,

$$
\begin{array}{ll} i \, d(B\circ Z)(E_{\alpha}, E_{- \alpha})  =&
-\frac{i}{2} B(Z, [E_{\alpha},
E_{-\alpha}])\\
=&-\frac{i}{2}B([Z,E_{\alpha}],E_{-\alpha})\\
= &-  \frac{i}{2}\alpha(Z)B(E_{\alpha}, E_{-\alpha})\\
=&- \frac{i\alpha(Z)}{<\alpha, \alpha>}\\
 = &-2i\frac{\alpha(Z)}{<\alpha, \alpha>} \omega_{\alpha} \wedge
\omega_{-\alpha}(E_{\alpha}, E_{-\alpha}).
\end{array}
$$

\begin{definition}\label{DefKosVect} 
The 1-form
$$\sigma = \sum_{\beta \in
R_{\mathfrak{m}}^+} \beta \in \gt^* \subset i\gc^*$$
 is called  the
Koszul form and  the dual vector $Z^{Kos}:= B^{-1} \circ \sigma$ is
called the Koszul vector.
\end{definition}

\begin{proposition}\label{propKoszul} \cite{A-P}  
The {\bf Koszul vector}
 $Z^{Kos}$
 defines  the invariant  K\"ahler-Einstein
 structure  $(\omega_{Z^{Kos}}, J(C))$ on $F=G/K$, where
$J(C)$ is the invariant complex structure  associated  with the
$T$-Weyl chamber $C$ which is defined by $\Pi_B$.
\end{proposition}

\smallskip

\noindent Let us conclude this section by recalling the flag manifolds of the classical groups:
(see, for example, \cite{A-P}, \cite{arv}):

\begin{enumerate}

\item[-] $SU(n)/S(U(n_1) \times \cdots \times U(n_s) \times U(1)^m)$

$n = n_1 + \cdots + n_s + m$, $s, m \geq 0$

\item[-] $SO(2n+1)/U(n_1) \times \cdots \times U(n_s) \times SO(2\ell+1) \times U(1)^m$

\item[-] $Sp(n)/U(n_1) \times \cdots \times U(n_s) \times Sp(\ell) \times U(1)^m$

\item[-] $SO(2n)/U(n_1) \times \cdots \times U(n_s) \times SO(2 \ell) \times U(1)^m$

$n = n_1 + \cdots + n_s + m + \ell$, $s, m, \ell \geq 0, \ell \neq 1$

\end{enumerate}

\subsection{Homogeneous  CR manifolds  }
We recall  that  a {\it CR  structure}  on a manifold $S$ is a pair
$(\mathcal{H}, J^{\mathcal{H}})$ where  $\mathcal{H}$ is a
codimension one  distribution and $J^{\mathcal{H}} \in
\Gamma(\mathrm{End}(\mathcal{H}))$ is a field of complex structures
in $\mathcal{H}$ such that the $\pm i$-eigendistributions
$\mathcal{H}^{\pm} \subset \mathcal{H}^{\mathbb{C}}$ are involutive
(i.e. closed  w.r.t. the Lie  bracket).\\

The   complex  structure   $J  $ on   a manifold  $M$  induces  a CR  structure  $(\ch,J^{\ch}) $ on
 any hypersurface $S  \subset M$  where
$\ch \subset TM$  is  the  maximal  $J$-invariant  distribution  and  $J^{\ch}= J|_{\ch}$.\\

\subsubsection{Ordinary   homogeneous  manifolds  and projectable  CR structures}

\begin{definition} 
A  homogeneous manifold $S=G/L$ of a  compact semisimple Lie
group $G$ is called  an  {\it ordinary  manifold}   if  the  normalizer $K:=
N_G(L)$ is the centralizer of a torus   and  contains $L$ as a
codimension one  subgroup.
\end{definition}

 Such a manifold is  the  total space  of  a  homogeneous  principal circle  bundle  $\pi: G/L \to  G/K$
 over  the flag manifold
$F =G/K$.\\

Let $S = G/L$ be an ordinary homogeneous manifold  and $K = N_G(L)$.
We define  the   {\it standard}  reductive  decomposition  as the
$B$-orthogonal  decomposition
 \be \label{reductivedecomposition}
\mathfrak{g} = \mathfrak{l}+ \mathbb{R}Z_F^0 + \mathfrak{m} \ee where
$\mathfrak{k} = \mathfrak{l}+ \mathbb{R}Z_F^0 $ and $B(Z_F^0,Z_F^0) = -1.$\\
We will call    $Z_F^0$ { \it  fundamental  vector }   and  we  will
identify it  with an invariant vector field on $S$ which is the
 fundamental vector  field of the
 principal bundle  $\pi$ which generates  a  commuting  with  $G$ action of   the circle   group $T^1$
 on $S$ , that  is  the  structure group  of  the principal bundle $\pi$.
The  dual 1-form  $\theta^0 : = B \circ Z_F^0$ is an invariant 1-form
on $S$ which  defines  a canonical invariant connection in $\pi$.

\noindent We denote by $Z^0$ the vector in $\mathfrak{t} = i Z(\mathfrak{k})$ such that $Z^0_F = i Z^0$ and, by a slight abuse of notation, we will call it also fundamental vector.

 The  following  lemma  shows  that, given a flag manifold  $F = G/K$, almost all  closed  codimension one subgroups $L$ of $K$ define  an ordinary manifold  $S = G/L$.

\begin{lemma}  \cite{A-S}, \cite{P-S}  
Let  $F = G/K$  be  a  flag manifold  and  $L$   a codimension one
 closed (normal)  subgroup of $K$.
If  $G/L$ is not  an ordinary manifold, then   $L = C_G(A_1)$   is  the centralizer  of  the 3-dimensional regular
 subalgebra $A_1$ of $\gg^{\mathbb{C}}$,
 associated    with a long root   such  that $G/ N_G(A_1)$ is  the  Wolf  space
 (symmetric   quaternionic K\"ahler manifold)   or $G = G_2$  and $\ga_1$ is the 3-dimensional subalgebra
 associated  with a short root.
 \end{lemma}

Now  we state    some elementary properties  of  an ordinary
manifold  $S = G/L$.
\begin{lemma} (\cite{A-S}, \cite{P-S}) Let $S = G/L$ an ordinary manifold. Then,
\begin{enumerate}
\item[i)]  Any invariant  vector  field on $S$ is proportional to  $Z^0_F$.
 \item[ii)] The  only  invariant
codimension one  distribution in   $S= G/L$  is the distribution
$\mathcal{H}$ defined  by the $\Ad_K$-invariant subspace
 $\mathfrak{m}$. This  distribution  is   also $Z^0_F$-invariant.
\item[iii)]  There is  a  natural one-to-one  correspondence between
invariant   complex  structures $J^F$ on  the  flag manifold  $F =
G/K$, $\mathrm{Ad}_K$-invariant   complex  structures  $J^{\gm}$  on
$\gm$ ( which are integrable in the  sense that
  $\gh^{\mathbb{C}}  + \gm^{10}$ is a subalgebra,  where  $\gm^{10} \subset \gm^{\mathbb{C}}$ is
 the $i$-eigenspace of  $J^{\gm}$ )   and    invariant CR  structures  $(\mathcal{H}, J^{\mathcal{H}})$
  on $S$    which   are also $Z^0_F$-invariant.
\item[iv)]  Any   irreducible $K$-submodule   of $\gm = T_oF$ remains
irreducible as $L$-submodule. \end{enumerate} 
\end{lemma}

 Following  \cite{P-S}  we   give

\begin{definition}\label{DefstandardCR}  
An invariant CR  structure $(\ch, J^{\ch})$ on an ordinary manifold $S=G/L$ is called
{\it projectable}  if it is  $Z^0_F$-invariant  or, equivalently  if
 the  projection   $\pi :S=  G/L \to F = G/K$    is  a  holomorphic map   of the CR manifold $S$ onto the  flag manifold $F$ with  some  invariant complex structure $J^F$.\\
An ordinary manifold  $S = G/L$  with a projectable  CR  structure is called
 a {\it    homogeneous CR manifold of  standard  type } or {\it   standard CR manifold}.
\end{definition}

 The   following  lemma  gives a  sufficient  condition for an  invariant CR  structure on an ordinary manifold $S = G/L$ to be projectable.

\begin{lemma} (\cite{P-S}) If irreducible $\Ad_K$-submodules of $\mathfrak{m}^{\mathbb{C}}$
 are non-equivalent  as $\Ad_L$-submodules  then any invariant CR
 structure on an ordinary  manifold $S = G/L$  is projectable.
 \end{lemma}

 All  compact homogeneous Levi non-degenerate CR manifolds  with non projectable   CR structure
  have been classified by  \cite{A-S}. More precisely, they prove

 \begin{theorem} \cite{A-S}
 A  compact homogeneous Levi  non-degenerate non projectable CR manifold $(S=G/L, \mathcal{H},  J|_{\mathcal{H}})$
 is either the  sphere bundle of  a  compact rank one  symmetric  space (CROSS)  or
one  of the exceptional  homogeneous  CR  manifolds
$$ SU_n /T^1 \cdot SU_{n-2},\,\,\,  SU_p  \cdot  SU_q/T^1  \cdot  U_{p-2}  \cdot  U_{q-2},$$
$$ SU_n/T^1  \cdot  SU_2  \cdot  SU_2  \cdot  SU_{n-4},\,\,\,  SO_{10}/T^1  \cdot  SO_6, \,\,\,  E_6/T^1  \cdot  SO_8 $$
which admit a holomorphic  fibration over  flag manifolds with
fibers  $S^n,\, n = 2,3,5,7,9$ respectively. 
\end{theorem}

\section{Cohomogeneity one K\"ahler manifolds  of  standard type}

Let $(M , \o, J, g)$ be  a  \C1 K\"ahler  $G$-manifold,
where $G$ is  a  connected compact semisimple Lie group.
 Deleting   singular  orbits  (if  they exist)   or,  in the  case
$M/G \simeq S^1$, the regular  orbit $S_0 = \pi^{-1}(0)$   we   get an open dense   submanifold  $M_{reg}$  of  regular points which   we identify   as a $G$-manifold
 with   $M_{reg} = (0,d) \times  G/L$, $d=1$ or $\infty$.\\

We  may assume   that  the  induced  metric $g_{reg}:= g|_{M_{reg}}$
has  the form
 $$    g_{reg} =  dt^2 + g_t,\,\,\,  t \in (0,d), $$
 where   $   g_t :=  g|_{S_t}$  is  an invariant  metric  on the
 regular orbit  $S_t := \{  t\} \times G/L = G/L$. We  denote  by $T_t := \partial_t|_{S_t}$
the  unit  normal vector  field to the orbit $S_t$    which is tangent  to normal geodesics
 $\gamma(t) = (t, x_0)$. Then $JT_t$ is a
$G$-invariant tangent vector  field on  $S_t$ and $\theta_t:= g\circ
JT_t$ the
dual invariant 1-form.\\
  Any regular orbit $S_t$ carries    an invariant    CR  structure  $(\ch, J^{\ch})$ induced  by $J$, where

 $$\mathcal{H} := JT_t^\perp = \mathrm{ker}\, \theta_t  \subset TS_t  $$
 and  $J_t := J|_{\mathcal{H}}$.

\begin{definition}\label{defSTANDARD}  
We  say that a \C1 K\"ahler  or  complex  manifold  $M$  is of  standard  type  if   a regular orbit $S_t = G \gamma(t) = G/L$
with the induced CR structure  is a  CR manifold of   the  standard  type (Definition \ref{DefstandardCR} above)   and   singular orbits (if  exist)  are flag manifolds  with
induced complex    structure.
\end{definition}

\begin{remark}
The last  condition is  automatically  satisfied if   $M$ is   a K\"ahler manifold.
\end{remark}

\subsection{  Moment map of a  \C1 \Km}

Let   $(M,J, \omega, g)$ be a \C1 \Km  of  a  compact semisimple  Lie group $G$. Since   the   group $G$ is  semisimple and preserves   the K\"ahler
 form, the
($G$-equivariant) moment map
 $$\mu : M \to \gg^* \simeq \gg,\,\,   x \to \mu_x,\,\, \mu_x(X) =h_X(x) $$
where $h_X$  it  the hamiltonian of  the Killing  vector  field $\hat X$
generated  by   $X \in \gg$,  is defined.

\begin{lemma} \cite{P-S}  
The restriction $\mu_t : S_t = G/L \to F_t :=\mu(S_t) =G/K$ of
the  moment map to a regular orbit  is a $G$-equivariant principal
bundle with   the  structure  group $T^1 =  K/L$    over a flag manifold  $F = G/K$
 (i.e.  an adjoint  orbit of
 $G$). In particular, $L$ is a  codimension one  normal  subgroup
 of  the  group $K$,  which is   the  centralizer  of   a torus in
 $G$.
 \end{lemma}

\pf Since the  moment map $\mu_t$ is  equivariant, it maps $S_t =
G/L$ onto an adjoint  orbit $F = G/K$.

   Now  $\ker \o_t = \ker \o|_{TS_t} = \mathbb{R}JT_t$  and   $d(h_X) = \o\circ X$
   give
   $\ker d\mu_t = \mathbb{R}JT_t$. So $L$ is a  codimension one
   subgroup of $K$. In terms of the Lie  algebra, we  can write a
   $B$-orthogonal  decomposition of $\gg$

   $$   \gg = \gl + \mathbb{R}Z_F^0 + \gm,   \,\, \gk = \gl + \mathbb{R}Z_F^0  $$
   where   $Z^0_F$ is  an $\mathrm{Ad}_L$-invariant vector, such that  the  associated
   invariant vector  field  on   $S_t$    is  proportional to  $JT_t$.
   The corresponding 1-parameter  group $T^1 = \mathrm{exp}\mathbb{R}JT_t = \mathrm{exp}\mathbb{R}Z^0_F$
(which is  closed  in  $K$, since  the   subgroup $L$ is  closed, see \cite{A-S})
   commutes  with $G$ and  defines  the structure of $T^1$-principal
   bundle $\pi : S_t \to F = G/K_t$.
   \qed

\medskip

\begin{remark} 
The  subgroup  $L \subset  K = L \cdot T^1  = L^{con} \cdot T^1$ can be  non connected.
\end{remark}

\begin{corollary}  
A cohomogeneity one  K\"ahler manifold $M$ is  of the  standard
type if the moment map $\mu :S_t = G/L \to F = G/K$ maps  any
regular orbit  with the induced CR structure   holomorphically onto the
associated  flag manifold  $F = G/K = G/N_G(L)$ with a  fixed
invariant complex structure $J^F$.
\end{corollary}

 \section{ Standard invariant  K\"ahler  structures on \C1 manifolds}
Let  $(S=G/L, \ch, J^{\ch})$ be a standard CR  manifold  with
holomorphic  fibration  $\pi : S =G/L \to F =G/K$  over  a  flag
manifold   with   a  complex  structure $J^F$. Infinitesimally,  it
is described  by  a  standard decomposition
 $$\gg = (\gl + \bR Z^0_F )+ \gm  =  \gk + \gm$$
where $Z^0_F$ is  the  fundamental vector, $\gk = N_{\gg}(\gl)$,
together   with  an $\mathrm{Ad}_K$-invariant complex  structure
$J^{\gm}$ in $\gm$. Now  we describe invariant Hermitian  and
K\"ahler structures of standard  type  in a regular \C1 manifold
$M_{reg} = (0,d)\times G/L$ which induce  a given  CR structure
$(\ch, J^{\ch})$ on each regular orbit  $S_t = \{t\}\times G/L =
G/L$.
  We identify $Z^0_F$ with a $G$-invariant vector
field on $S_t = G/L$ and denote  by $\theta^0 := B \circ Z^0_F$   the
dual invariant  1-form   with kernel $\ch = \ker \theta^0$. Any
invariant vector field  (resp. 1-form)  is  conformal  to $Z^0_F$
(resp. $\theta^0$).

\subsection{Invariant   Hermitian   structures on   $M_{reg}$}
 We   fix  an invariant metric $g$  on  $M_{reg}= (0,d) \times G/L$  of  the  form
$$   g = dt^2 + g_t$$
where $g_t$ is  an invariant metric in $S_t = S = G/L$  such that the  CR  structure $J^{\ch}$
is  orthogonal.\\
  We  describe all  $g$-orthogonal  projectable invariant  complex  structures $J$ on $M_{reg}$ which
  project onto a  given invariant complex  structure $J^F$ of   the  flag manifold $F = G/K$,
  or, equivalently, are extensions  of    the associated  CR   structure  $J^{\ch}$.
     Then    $(g, J)$  is   an invariant  Hermitian structure  in $M_{reg}$.\\

Since  the  invariant 1-forms   $dt , \theta^0$   are orthogonal to
each other  and the  distribution $\ch = \mathrm{ker}\{dt, \theta^0
\}$, any   extension of  $J^{\ch}$ to  an  orthogonal invariant
  almost complex  structure on $M_{reg}$ can be written as
\be \label{complexstructure}
  \begin{array}{lllll}
 J^* dt &:=& dt \circ J= a(t) \theta^0,\,\, & J^*
\theta^0 =& -\frac{1}{a(t)}dt\\
J T_t &=& \frac{1}{a(t)}Z^0_F,\,\, &JZ^0_F=& -a(t) T_t
  \end{array}
\ee
  where  $a: (0, d) \to \bR$  is  a non vanishing  function.
	
\begin{proposition}
Any extension of  the  CR  structure  $J^{\ch}$   to   an invariant
orthogonal   (integrable) complex  structure $J$  on $M_{reg}$ is
given by (\ref{complexstructure}).
\end{proposition}

\pf We have    to check that   the  almost complex  structure $J$ is
integrable.
 Since $ J|_{\ch} = J^{\ch}$  is   an integrable CR  structure  in $ S$,
  it is  sufficient  to check that   the differential of the 1-form
  $dt - iJ^* dt = dt + ia   \theta^0  \in \Omega^{1,0}(M) $ belongs  to  the
    space $\Omega^{2,0}(M) + \Omega^{1,1}(M)$ of   forms of  type $(2,0)$  and $(1,1)$.
Since $d \theta^0  = \omega^0 \in \Omega^{1,1}(M)$ , we  have
$$
\begin{array}{ll}   d( dt + ia    \theta^0 )    &= i \dot a dt \wedge \theta^0  +  ia d \theta^0 \\
 &= (dt +i a    \theta^0) \wedge i \dot a  \theta^0 +   ia d \theta^0 \in \Omega^{2,0}(M) + \Omega^{1,1}(M).
\end{array}
$$

\qed

\subsection{Invariant   K\"ahler   structures on   $M_{reg}$}

We describe all    standard   invariant  K\"ahler   structures  $(\omega_{reg}, J_{reg}, g_{reg})$ on
 $M_{reg}= (0,d) \times G/L$  which  induce the projectable CR structure   $(\ch, J^{\ch})$
on $S_t = \{t\}  \times  S$.

Recall  that
vectors $Z \in \gt = i Z(\gk)$ correspond to invariant closed 2-forms $\o_Z$ on
$F= G/K$, whose value   at $o = eK$  is given by

$$ (\o_Z)_o (X,Y)  =i \ d (B \circ Z)(X,Y), \ X, Y \in \gm.$$
  We  will denote  the pull back of $\omega_Z$   to $M_{reg}$ by the  same letter  $\o_Z$.
In  particular, the  vector  $Z^0 = -iZ^0_F$   defines  the  invariant
2-form $\o^0:= \o_{Z^0} = i \ dB\circ Z^0 = dB \circ Z^0_F = d \theta^0$ on $S$ which
is the curvature of  the principal  connection $\theta^0$ in the
principal bundle $ \pi : S = G/L \to F = G/K $. We denote by $C(J^F)
\subset  \gt$    the $T$-Weyl chamber which corresponds to the
complex  structure $J^F$.

The  following proposition shows  that a standard invariant
K\"ahler structure in $M_{reg}$  is  defined   by a parametrized
open interval $(Z_0 Z_d)$ in the $T$-Weyl chamber $C(J^F)$, parallel
to $Z^0$.

 \begin{proposition}\label{KaehlerstructonM{reg}} 
A  standard invariant K\"ahler
structure   $(\o_{reg}, J_{reg},  g_{reg})$ on $M_{reg} = (0,d)
\times G/L$    which induces an
   invariant  CR structure $(\ch, J^{\ch})$ on regular orbits  is defined  by  a
parametrized open interval $Z_t= Z_0 + f(t)Z^0 \subset C(J^F)$ where
$f: (0,d) \to \mathbb{R},\, \lim_{t \rightarrow 0 }f(t) =0$ is a   smooth function  with
$a(t)=\dot f(t)>0$. More precisely,  define
$$\theta_t := iB\circ Z_t = iB\circ Z_0 +f(t) iB\circ Z^0 = \theta_0 + f(t) \theta^0,$$
$$ \omega_0 := d\theta_0 = \omega_{Z_0}, \,  \omega^0:= d\theta^0 = \omega_{Z^0}.$$
 Then the  K\"ahler  structure is  given by
 $$  \omega_{reg} = d(\theta_0 + f(t) \theta^0) = \dot f dt \wedge \theta^0 + \omega_0 + f(t) \omega^0, $$
$$ g_{reg} = dt^2 + (\dot f \theta^0)^2 + \pi^* g_0 + f(t) \pi^*g^0,
$$
$$ J_{reg}: dt \to - \dot f \theta^0,  \theta^0 \to \frac{1}{\dot f} dt, \, J|_{\ch} = J^{\ch}. $$
Here $ g^0 = -\omega^0 \circ J^F, \,  g_0 = -\omega_0 \circ J^F$ are   symmetric  bilinear forms on  $F$.\\
The pair $(\omega_t:= \omega_0 + f(t)\omega^0, J^F)$  defines an
invariant K\"ahler structure on $F$ for  $t \in (0,d)$ .
\end{proposition}

 \pf   Using (\ref{complexstructure}), we have
$$\omega_{reg}(T_t, JT_t) =1 =\omega_{reg}(T_t, \frac{1}{a}Z^0_F).$$

\noindent This  shows   that the K\"ahler   form $\o_{reg}$ on $M_{reg}$ can
be written as
$$   \o_{reg} = dt\wedge a(t)\theta^0 + \omega^S_t $$
where $\omega^S_t := \omega|_{S_t}  $  is  a closed  invariant
2-form on $S_t =S$  with  kernel $\bR Z^0_F$.

\noindent The  form $\omega_t^S$ is the pull back  of   an invariant
 symplectic  form  $\omega_{Z_t} = i d(B\circ Z_t)$  associated  with
 a  vector  $Z_t \in C(J^F)$ (where $C(J^F)$ is the  Weyl chamber which corresponds to
 the  complex  structure $J^F$). It is  sufficient  to check  that $\omega_t^S$ is invariant with
respect to  the  fundamental  vector  field $Z^0_F$. We have
$$ \cl_{Z^0_F} \omega_t^S = i_{Z^0_F}d \omega_t^S + di_{Z^0_F}\omega_t^S =0$$
since  $\omega_t^S$  is  closed  and  has kernel $\mathbb{R}Z^0_F$.

\noindent Now  the  condition that $\omega_{reg}$ is  closed  can be  written
as
$$dt\wedge a \omega^0 + dt \wedge \dot \omega_t^S =0$$
or

$$    -a(t) \omega^0 + \dot{\omega}^S_t =0 = \omega_{a Z^0 + \dot Z_t} =0$$
or  $\dot Z_t + aZ^0 =0.$
 This  implies  that   the  curve $Z_t $  is an (open) interval
 (maybe, a ray)   $(Z_0 Z_d) \subset \bar C(J^F) $ with parametrization
$Z_t = Z_0 + f(t) Z^0, \, t \in (0,d)$, with $a(t) = \dot f(t)
>0$. We may assume that   $f(0) =0$. We have  proved  that
$$\o_{reg} = \dot f(t) dt \wedge \theta^0   +  \omega_{Z_t}  = \dot f(t) dt \wedge \theta^0 +    \o_0 + f(t) \o^0 .$$

Now  we  easily calculate  the  metric $g_{reg}$ as a composition of
$\omega_{reg}$ and $J$. \qed \\

Since an interval in $C(J^F)$ is not  diffeomorphic  to  a circle,
we get

\begin{corollary}
There is no \C1  \Km  of  standard  type  with the orbit
space $S^1$.
\end{corollary}

 \begin{corollary}\label{corollarya(t)andb(t)}
 \begin{enumerate}
   \item[i)]  $a(t):= \dot f(t) = \omega_{reg}(T_t, Z^0_F) = g_{reg}(JT_t,  Z^0_F),\,$\\
   $ a(t)^2 = g_{reg}((Z^0_F)_{\gamma(t)}, (Z^0_F)_{\gamma(t)} ). $
   \item [ii)] For any $X \in \mathfrak{g}$ the square norm of the
   Killing field $\hat X$ along a normal geodesic $\gamma(t)$
 is given by
    $$  b_X(t)= |\hat X|^2_{\gamma(t)} = \omega_0(X,JX) + f(t) \omega^0(X,JX) =  g_0(X,X)+ f(t) g^0(X,X).$$
   \item[iii)] If  $\omega^0(X,JX) \neq 0$,  then $\dot b_X(t) = a(t) \omega^0(X,JX) = a(t)g^0(X,X)\neq
   0$ and  the  function $b_X(t)$
   has no critical points for  $t \in (0,d)$. It is true if $0 \neq X \in
   \mathfrak{m}$.
 \end{enumerate}
\end{corollary}

\medskip
\subsection{Basic properties  of standard  K\"ahler \C1 manifolds  with  singular orbits}
As  an application of previous results, we  prove  two  basic
properties  of a  standard  K\"ahler \C1 manifold  with one or two
singular orbits (see also \cite{P-S}).

\begin{proposition}\label{Prop  singular orbit }  
Let $M$  be a standard  \C1 K\"ahler manifold  with
  the orbit  space $M/G = [0,1)$ .
Then    the  singular  orbit
$S_0 = G \gamma(o) = G/H_0$   is a  complex  submanifold, hence  a  K\"ahler   flag  manifold and
 $H \supset K = N_G(L)$.
\end{proposition}

\pf The value $\hat Z^0_p$  of the   Killing  vector  field $\hat Z^0$
generated   by the  fundamental vector $Z^0_F \in \gg$   at the point
$p:=\gamma(0) \in S_0$ is zero, since in  the opposite case we get
two $\Ad_L$-invariant elements $\hat Z^0_p, J\hat Z^0_p $ in the
$B$-orthogonal  complement  $\mathfrak{l}^\perp  = \bR Z^0_F + \gm
\subset \gg$. This proves that $H
\supset K$   since  $K =  N_G(L)$ is  a  connected  subgroup. \\
 For a  subspace
 $\mathfrak{n} \subset \mathfrak{g}$ we denote by $\hat{\mathfrak{n}}_{\g(t)}$ the subspace of
$T_{\gamma(t)}M$   spanned by  the values   of  the Killing  vectors  ${\hat X_{\gamma(t)}}, \, X \in \gn$.
Since  $\omega (\dot\gamma(t), \hat \gm_{\g(t)})  = 0$ for  $ t \neq 0$  it is true  also for  $t=0$. But
$$    T_{\gamma} S_0 = \hat \gg_0  = \hat \gm_{\g(0)}$$
since  $(\hat \gl +\bR \hat Z_0)_0  =0$. So   for   any normal  geodesic  $\gamma(t)$   with $\gamma(0) =p$, we  get
$$ \omega (\dot \gamma(0), T_p S_0) =g(\dot \gamma(0), JT_p S_0)=   \omega(\dot \gamma(0),  \hat \gm_{\g(0)})  = g(\dot \gamma(0), J \hat \gm_{\g(0)}) =0.$$
Since  vectors  $\dot \gamma(0)$  span $T^{\perp}_p S_0$,   we get
$$  g(T^{\perp}_p S_0, J  T_p S_0) =0$$
which shows  that  $S_0$ is  a  complex  submanifold. \qed\\

 As a corollary, we get the

\begin{proposition} 
If $(M,\omega, J, g)$ is a standard  (compact) \C1 \Km with two
singular orbits  $S_{\epsilon} = G \gamma(\epsilon) =
G/H_{\epsilon}, \, \epsilon = 0,1$, then  the   singular orbits are
complex submanifolds  and   $ K:=N_G(L) = H_0 \cap H_1$.
\end{proposition}

\pf It remains  only to check that   $H_0 \cap H_1  = K $. Since the
metric $g$ is  complete,  for  any normal geodesic  $\gamma(t),\, t
\in \bR$  through  a point $p \in M$  we have
$$ \omega(\dot\gamma(t) ,\hat{\mathfrak{m}_t}) =0.$$
 Also, since $\hat
Z^0_p=0$, we get
$$\omega(\dot \gamma(0), \hat Z^0_{p}) = 0.$$
 Since $T_pS_0 =
\hat{\mathfrak{m}}_0 $, this means  that $\omega(T^{\perp}_p S_0,
T_pS_0) = g(T^{\perp}_p S_0, JT_pS_0))=0$ that is $S_0$ is  a
complex  submanifold and $H_0$ and similarly $H_1$ are the
centralizers of some torus  in $G$. Then $H_0 \cap H_1 \supset K$ is
also  the centralizer of  a torus and hence, a connected  subgroup.
If $K \neq H_0 \cap H_1$ , then there  is non-zero vector $X \in
\mathfrak{m}\cap \mathfrak{h_0}\cap \mathfrak{h}_1$ and the
associated   function $b_X(t):= |\hat X_{\gamma(t)}|^2$ (the square
norm of  the Killing  field $\hat X$ along  a  normal geodesic)
vanishes  at the points $t =0,1$. Hence it has  a critical point in
$(0,d)$, which contradicts Corollary \ref{corollarya(t)andb(t)}. \qed

\section{Einstein equation and  K\"ahler-Einstein structure on
$M_{reg}$} \label{einstequation}

To calculate   the  Ricci form  of an invariant K\"ahler metric  on
$M_{reg} = (0,d) \times G/L$  we construct  local  holomorphic
coordinates
 $z^0, z^1, \dots, z^n$  in  $M_{reg}$  which  are   extension  of local  holomorphic  coordinates
  $z^1, \dots, z^n$ of  the associated  flag manifold  $F = G/K$.

\subsection{Holomorphic  coordinates  and K\"ahler potential of   $M_{reg}$}

Let $z_1, \dots, z_n, t, \phi$ be local coordinates on $M_{reg}$, where $\phi$ is a local coordinate on the torus $T^1 = \{ e^{i \phi} \}$ such that $K = L \cdot T^1$. In these coordinates, let

\begin{equation}\label{thetazero}
\theta^0 = c \ d \phi + \Phi(z, \bar z) = c \ d \phi + i \sum_{j=1}^n(F_j d z_j - \bar F_j d \bar z_j)
\end{equation}

\noindent and let $\Psi$ be a solution to the system of partial differential equations

\begin{equation}\label{diffeq}
\frac{\partial \Psi}{\partial z_j} = F_j, \ \ j = 1, \dots, n
\end{equation}

\noindent (it is easily checked that $d \theta^0 = \omega^0$ implies $\frac{\partial F_k}{\partial z_j} = \frac{\partial F_j}{\partial z_k}$).

\noindent Then, we claim that

\begin{equation}\label{z0}
z_0 = i c \phi + \int \frac{1}{f'} dt + \Psi
\end{equation}

\noindent is a new holomorphic coordinate on $M_{reg}$. Indeed, this is equivalent to say that $\bar \partial z_0 = (d - i d^c)(z_0)=0$, where $d = \partial + \bar \partial$ and $d^c = i(\bar \partial - \partial) = J^{-1} d J$. We have

$$dz_0 =  i c \ d \phi + \frac{1}{f'} dt + \sum_{j=1}^n \frac{\partial \Psi}{\partial z_j} dz_j + \frac{\partial \Psi}{\partial \bar z_j} d \bar z_j=$$
$$= i c \ d \phi + \frac{1}{f'} dt + \sum_{j=1}^n F_j dz_j + \bar F_j d \bar z_j $$

\noindent and

$$d^c z_0 = J^{-1}( i c \ d \phi + \frac{1}{f'} dt + \sum_{j=1}^n F_j dz_j + \bar F_j d \bar z_j) =$$
$$= -i c J d \phi - \frac{1}{f'} J dt - i(\sum_{j=1}^n F_j dz_j - \bar F_j d \bar z_j) =$$

\noindent (by $J dt = - f' \theta^0$ and $J \theta^0 = \frac{1}{f'} dt$)

$$= -i \frac{1}{f'} dt - i (\sum_{j=1}^n F_j dz_j + \bar F_j d \bar z_j) + \theta^0 - i(\sum_{j=1}^n F_j dz_j - \bar F_j d \bar z_j) =$$
$$= -i \frac{1}{f'} dt + c \ d \phi - i(\sum_{j=1}^n F_j dz_j + \bar F_j d \bar z_j) $$

\noindent so $dz_0 = i d^c z_0$ as required.

\subsection{The Ricci form  and   the Einstein equation  for  an invariant K\"ahler metric  on $M_{reg}$}\label{einstequation2}

Now, we calculate  the Ricci form  of  an  invariant K\"ahler
structure $( \omega_{reg}, J^{reg}, g^{reg})$  on $M_{reg} = (0,d)
\times G/L$  associated with a parametrized interval
  $ Z_t = Z_0 + f(t)Z^0$ in the $T$-Weyl chamber $C(J^F)$ corresponding to the complex structure $J^F$ of the flag manifold $F = G/K$.

\noindent By Proposition \ref{KaehlerstructonM{reg}}, the  K\"ahler  form is   given by

\begin{equation}\label{omegaregbyProp4}
 \omega_{reg} = \dot f(t) dt\wedge \theta^0 + (\omega_0 + f(t) \omega^0) = \dot f(t) dt\wedge \theta^0  + \omega_{Z_t}
\end{equation}

where $\omega_{Z_t} = d(B\circ Z_t)  $, so that

$$\omega_{reg}^{n+1} = (n+1) df \wedge \theta^0 \wedge
(\omega_{Z_t})^n  , $$

$$
\begin{array}{ll}
 (n+1)(\omega_{Z_t})^n =& (n+1)(-2)^n(\sum_{\alpha \in R_{\mathfrak{m}}^+} \frac{\alpha(Z_t)}{(\alpha, \alpha)}
 \omega_{\alpha} \wedge \omega_{-\alpha})^n =  h(f) \mathrm{vol}^F,
\end{array}
$$
where $R_{\mathfrak{m}}^+$ is the set of positive black roots,
 $$\mathrm{vol}^F := (n+1)! (-2)^n
\prod_{\alpha \in R_{\mathfrak{m}}^+} \frac{1}{(\alpha,
\alpha)}\omega_{\alpha} \wedge \omega_{-\alpha} $$ is  a volume form
on $F = G/K$ and
$$ h(f) =
\prod_{\alpha \in R_{\mathfrak{m}}^+}\alpha(Z_t) = \prod_{\alpha \in
R_{\mathfrak{m}}^+}(\alpha_0 + f \alpha^0),\,\, \, \alpha_0 :=
\alpha(Z_0),\, \alpha^0 := \alpha(Z^0).$$
  So we can finally write

\begin{equation}\label{omegavolume}
\omega_{reg}^{n+1} = \dot f h(f)  dt \wedge \theta^0 \wedge \mathrm{vol}^F.
\end{equation}

\begin{lemma} 
Let $z^1, \dots, z^n$ be  holomorphic  coordinates in $F =
G/K$ and $z^0, z^1, \cdots, z^n$ their extension to  holomorphic
coordinates in $(0,d) \times G/L$, being $z_0$ defined by above formula (\ref{z0}).  Then  the  density $\mu$
associated with the   volume  form
$$\omega_{reg}^{n+1} = \mu dz^0 \wedge \cdots \wedge dz^{n} \wedge d \bar z^0 \wedge \cdots  \wedge d \bar z^{n}$$
 is  given by

$$ \mu = \frac{i}{2} \dot f^2 h(f) \mu^F $$
 where $\mu^F$ is the  density associated   with the  volume  form
 $\mathrm{vol}^F$.
\end{lemma}

 \pf
 By substituting (\ref{thetazero}) in (\ref{omegavolume}) we have

\begin{equation}\label{densityomegareg}
\omega_{reg}^{n+1} = \dot f h(f) \mu^F dt \wedge d \phi \wedge dz^1 \wedge \cdots \wedge dz^{n} \wedge d \bar z^1 \wedge \cdots  \wedge d \bar z^{n}
\end{equation}

By (\ref{z0}) we have $d z_0 = i c \ d \phi + \frac{1}{f'} dt +$ (terms containing $dz_j$, $d \bar z_j)$, $d \bar z_0 = -i c \  d \phi + \frac{1}{f'} dt +$ (terms containing $dz_j$, $d \bar z_j)$ and $ -2i c \frac{1}{f'} dt \wedge d \phi =  dz_0 \wedge d \bar z_0 +$ (terms containing $dz_i$, $d \bar z_i$). By replacing this in (\ref{densityomegareg}) we conclude the proof.

 \qed

Now  we calculate  the Ricci  form  $\rho$  of   the  invariant
K\"ahler   structure using the very well-known formula 
$$\rho = -i \partial \bar \partial \log \mu = -\frac12   d d^c \log
\mu  = -\frac12   d d^c \log (\dot f^2 h(f)) - \frac12   d d^c \log
\mu^F.$$
 The second  term is  the Ricci  form $\rho^F = d \sigma$ of (any) invariant
 K\"ahler (with  respect to  the  complex structure $J^F$)  metric.
 We calculate the first term. By Proposition \ref{KaehlerstructonM{reg}}
 $$ d^c\log (\dot f^2  h) = J^{-1}  d\log (\dot f^2 h)= \frac{d}{dt}\log (\dot f^2 h) J^{-1} dt = \dot f \frac{d}{dt}\log (\dot f^2 h)\theta^0 ,$$
$$ d d^c \log (\dot f^2  h)= \frac{d}{dt}( \dot f \frac{d}{dt}\log (\dot f^2 h)) dt \wedge
\theta^0 + \dot f \frac{d}{dt}\log (\dot f^2 h)\omega^0
 $$
\noindent and
$$\rho = -\frac12 \frac{d}{dt}( \dot f \frac{d}{dt}\log (\dot f^2 h)) dt \wedge
\theta^0 - \frac12 \dot f \frac{d}{dt}\log (\dot f^2 h)\omega^0 +
\rho^F. $$ The Einstein equation $\rho = \lambda \omega_{reg} \equiv
\lambda \dot f dt\wedge \theta^0 + \lambda \omega_0 + \lambda f
\omega^0$ can be rewritten as

\begin{equation}\label{einsteinfirstequation}
\lambda f + \frac12 \dot f \frac{d}{dt} \log(\dot f^2 h)= c, \, c= const
\end{equation}
\begin{equation}\label{einsteinsecondequation}
\rho^F = \lambda \omega_0 + c\omega^0
\end{equation}
Since by Proposition \ref{propKoszul} we have  $\rho^F = d (B \circ Z^{Kos})$, where
$Z^{Kos}$ is the Koszul vector, (\ref{einsteinsecondequation}) can
be rewritten as
 \be \label{Einstein1}
  Z^{Kos} = \lambda Z_0 + cZ^0.
  \ee
  Now  we write (\ref{einsteinfirstequation}) in more explicit  form.
 Since $ \dot f \frac{d}{dt} \log(\dot f^2) = 2 \ddot f $ and
 $$ \frac{d}{dt} \log( h(f))= \frac{d}{dt} \sum_{\alpha \in
 R_{\mathfrak{m}}^+}\log (\alpha_0 +f \alpha^0)=\sum_{\alpha \in
 R_{\mathfrak{m}}^+} \frac{\dot f \alpha^0}{\alpha_0 + f
 \alpha^0}=:\dot f A(f)
$$
where $A(f) = \sum_{\alpha \in
 R_{\mathfrak{m}}^+} \frac{\alpha^0}{\alpha_0 +f \alpha^0}$, we get
 \be \label{Einstein2}
\ddot f + \frac12 A(f)\dot f^2 + \lambda f  =c \ee
 We proved the following

 \begin{proposition} \cite{P-S} 
The K\"ahler metric on $(0,d) \times G/L$ defined  by a
 parametrized interval $ Z_t = Z_0 + f(t)Z^0 \subset C, \, t \in
 (0,d)$ of  a Weyl chamber $C$ is a K\"ahler-Einstein metric if and only if the
  function $f(t)$ satisfies  the  equation (\ref{Einstein2}) and
  the  relation (\ref{Einstein1}) holds.
 \end{proposition}

\subsection{An example: the $n$-dimensional projective space}\label{sectionexample}

Let $M$ be the projective space $\bC\bP^n = \{ [x_0, x_1, \dots x_n] \}$ with affine coordinates $w_k = \frac{x_k}{x_0}$, $k = 1, \dots, n$, endowed with the Fubini-Study form $\omega_{\bC\bP^n} = i \partial \bar \partial \log (1 + |w_1|^2 + \cdots + |w_n|^2)$. Let us denote $x = (x_1, \dots, x_n)$ and let us consider the action of $G = SU(n)$ on $M$ given by $A[x_0, x] = [x_0, Ax]$.

\noindent It is easily checked that this is a cohomogeneity one action, with singular orbits $S_0 = G([x_0, 0]) = {[x_0, 0]}$, $S_1 = G([0, x]) = \bC\bP^{n-1}$, and regular orbits diffeomorphic to $S^{2n-1}$ (obtained when $x_0 \neq 0$ and $x \neq 0$).

\noindent For every regular point $[1, w_1, \dots, w_n]$, its orbit can be identified with the sphere of radius $r = \sqrt{|w_1|^2+ \cdots + |w_n|^2}$ in $\bC^n$, so that the singular orbit $S_0$ (resp. $S_1$) is obtained when the radius tends to 0 (resp. to infinity).

\noindent The flag $F = G/K$ associated with the regular orbits is $SU(n)/U(n-1) = S^{2n-1}/\{e^{i \phi}\} = \bC\bP^{n-1}$. Let $z_k = \frac{w_k}{w_1}$, $k=2, \dots, n$, be the affine coordinate on $F$. Then we can take $r, \phi, z_2, \dots, z_n$ as local coordinates on a dense open subset of the union of regular orbits $M_{reg}$. More precisely, set

$$(w_1, \dots, w_n) = \left( r \frac{w_1}{\sqrt{|w_1|^2+ \cdots + |w_n|^2}}, \dots,  r \frac{w_n}{\sqrt{|w_1|^2+ \cdots + |w_n|^2}} \right) = $$
\begin{equation}\label{coordinates}
= \left( \frac{r e^{i \phi}}{\sqrt{1+|z|^2 }}, \frac{ r e^{i \phi} z_2}{\sqrt{1+|z|^2 }}, \dots, \frac{ r e^{i \phi} z_n}{\sqrt{1+|z|^2 }} \right).
\end{equation}

\noindent where we are setting $|z|^2 := |z_2|^2 + \cdots + |z_n|^2$.

\noindent Using the change of coordinates (\ref{coordinates}), after a long but straight computation we can see that the restriction $\omega_{reg}$ of the Fubini-Study form $\omega_{\bC\bP^n}$ to $M_{reg}$ is

\begin{equation}\label{omegaregrzphi}
\frac{2r}{(1+r^2)^2} dr \wedge \left( d \phi + \frac{i}{2(1 + |z|^2)} \sum_{j=2}^n (z_j d \bar z_j - \bar z_j dz_j) \right) + \frac{r^2}{1+r^2} \tilde \omega
\end{equation}

\noindent where $\tilde \omega = i \partial \bar \partial \log (1 + |z|^2)$ is the Fubini-Study form on $\bC\bP^{n-1}$.

\noindent In order to find the relation between the parameter $r$ and the parameter $t$ used in the above sections, let us recall that $\frac{\partial}{\partial t}$ is a unit field on $M_{reg}$ normal to each regular orbit. On the one hand, $\frac{\partial}{\partial r}$ is normal to the regular orbits. On the other hand, a straight calculation shows that the coefficient of $dr^2$ in the espression of the metric $g_{reg}$ is $\frac{2}{(1+r^2)^2}$, so we must have $dt = \frac{\sqrt 2}{1+r^2} dr$ and, integrating, $r = \tan \left( \frac{t}{\sqrt 2} \right)$. By replacing this into (\ref{omegaregrzphi}) we find

\begin{equation}\label{omegaregrzphi2}
\frac{\sqrt 2}{2}\sin \left( \frac{t}{\sqrt 2} \right) dt \wedge \left( d \phi + \frac{i}{2(1 + |z|^2)} \sum_{j=2}^n (z_j d \bar z_j - \bar z_j dz_j)  \right) + \sin^2 \left( \frac{t}{\sqrt 2} \right) \ \tilde \omega
\end{equation}

\noindent This is exactly formula (\ref{omegaregbyProp4}) with $f(t) = \frac{\sqrt{n-1}}{n \sqrt 2} \sin^2 \left( \frac{t}{\sqrt 2} \right)$, $\omega^0 = n \sqrt{\frac{2}{n-1}} \tilde \omega$ and $\theta^0 = n \sqrt{\frac{2}{(n-1)}} \left( d \phi + \frac{i}{2(1 + |z|^2)} \sum_{j=2}^n (z_j d \bar z_j - \bar z_j dz_j) \right)$. The correct normalization for $\theta^0$ is obtained as follows. On the Lie algebra $su(n)$ of $G = SU(n)$, $\theta^0$ is defined as $B(Z^0_F, \cdot)$, where $B(X, Y) = 2n \tr(XY)$ is the Killing-Cartan form and $Z^0_F \in su(n)$ is such that $B(Z^0_F, Z^0_F) = -1$. This last condition easily implies that $Z^0_F = \frac{i}{n \sqrt{2(n-1)}}\diag((n-1), -1, \dots, -1)$. Then, for every $X_I \in T_I (SU(n))$ having $i\alpha_1, \dots, i \alpha_n$ ($\alpha_k \in \bR$) on the diagonal, we have
$$\theta^0(X_I) = B(Z^0_F, X_I) = -\frac{2n}{{n \sqrt{2(n-1)}}} ((n-1) \alpha_1 - \alpha_2 - \cdots - \alpha_n) =$$
$$= -\sqrt{\frac{2}{(n-1)}}n \alpha_1$$

\noindent and for general $X_g \in T_g (SU(n))$, $g = (a_{ij})$,

$$\theta^0(X_g) = \theta^0(dg^{-1} X_g) = -n\sqrt{\frac{2}{(n-1)}}(\bar a_{11} a_{11}' + \cdots + \bar a_{n1} a_{n1}' ).$$

\noindent Hence $\theta^0 = -\sqrt{\frac{2}{(n-1)}}(\bar a_{11} d a_{11} + \cdots + \bar a_{n1} d a_{n1} )$. Now, if we take into account that $\bC\bP^{n-1} = SU(n)/U(n-1)$ via the action $\bC\bP^{n-1} = \{ A [1, 0, \dots, 0]\} = \{ [a_{11}, \dots, a_{n1}] \}$ and accordingly to (\ref{coordinates}) we replace
$a_{11} = \frac{e^{i \phi}}{\sqrt{1+|z|^2 }}, a_{21} = \frac{ e^{i \phi} z_2}{\sqrt{1+|z|^2 }}, \dots, a_{n1} = \frac{ e^{i \phi} z_n}{\sqrt{1+|z|^2 }}$, a straight calculation yields exactly the above expression for $\theta^0$.

\noindent Finally, let us verify that $f(t) = \frac{\sqrt{n-1}}{n \sqrt 2} \sin^2 \left( \frac{t}{\sqrt 2} \right)$ satisfies the Einstein equation (\ref{Einstein2}). In this case, since $F = \bC\bP^{n-1}$, the set $R_{\mathfrak{m}}^+$ contains the roots $\varepsilon_1 - \varepsilon_2, \dots, \varepsilon_1 - \varepsilon_n$. Moreover, since the singular orbit $S^0$ is a fixed point, we must have $Z_0 = 0$ and then $\alpha_0 = \alpha(Z_0)=0$. Hence $A(f) = \sum_{\alpha \in R_{\mathfrak{m}}^+} \frac{\alpha^0}{\alpha_0 + f \alpha^0} = \frac{n-1}{f}$ and (\ref{Einstein2}) reduces to

$$\ddot f + \frac{n-1}{2} \frac{\dot f^2}{f} + \lambda f = c .$$

\noindent One immediately sees that $f(t) = \frac{\sqrt{n-1}}{n \sqrt 2} \sin^2 \left( \frac{t}{\sqrt 2} \right)$ is a solution to this equation for $c = \frac{\sqrt{n-1}}{\sqrt 2}$ and $\lambda = n+1$, which is exactly the value of the Einstein constant for $\omega_{\bC\bP^n} = i \log (1 + |w_1|^2 + \cdots + |w_n|^2)$.

\section{Standard cohomogeneity one  K\"ahler manifolds    with one  singular orbit}\label{Section C1standard}

In the   sequel  we   will   consider    only     standard \C1  $G$-manifolds  $M$ (Definition \ref{defSTANDARD})  with one (complex) singular orbit  $(S_0= G/H, J^S)$. For  brevity,  we will call these manifolds just {\it standard C1 manifolds}.

\smallskip

\noindent First of  all, we prove  that any such manifold is the  total
space  of the homogeneous  vector  bundle $M_{\vp} = G\times_H
V_{\vp}\to S_0 =G/H$  over   a flag manifold  $S_0 = G/H$ defined by
 a  representation $\varphi : H \to  GL(V_{\vp})$  with    $\vp(H) \simeq U_m$. Then we  give  a description of  these manifolds in terms of painted
  Dynkin  diagrams and  determine   the invariant K\"ahler
metrics  on them.

\subsection{Reduction to  admissible vector bundles}

Let  $(S_0 =G/H,J^S)$ be  a   flag manifold  with  invariant  complex  structure.

\begin{definition} 
A complex linear  representation  $\vp : H \to GL(V_{\vp}) =
GL_m(\mathbb{C})$  and the  associated  homogeneous  vector  bundle
$M_{\vp} = G \times_H V_{\vp} $   are called {\it admissible} if
$\vp(H) \simeq U_m$.
\end{definition}

\noindent Note  that  an  admissible  representation is defined  by  a normal
subgroup  $N = \mathrm{ker}(\vp)$ such that  $H/N \simeq U_m$.
Such a subgroup $N$   is also called an {\it admissible  subgroup}.\\

\noindent The main result of this section is the following

\begin{theorem}\label{C1  complex manifolds} 
Let  $M$ be a standard  complex  C1 manifold. Then $M$ is the total space $M_{\vp}$ of an admissible vector bundle $ M_{\vp} = G \times_H V_{\vp} \rightarrow S_0$ over the singular orbit $(S_0 = G/H, J^S)$.
\end{theorem}


\pf  Any \C1  manifold $M$  with one singular orbit can be identified  with   the
homogeneous vector  bundle $M = G \times_H V$  over
    the singular orbit $S_0 = G/H$ associated to some sphere transitive representation $\vp : H \to GL(V_{\vp})$,
     where  $V_{\vp}$  is  the normal space  of  the singular orbit $S$ (Proposition \ref{RiemC1manifolds}). Since, by assumption,  the orbit
     $S_0$  is a complex submanifold, the linear group $\vp(H)$ preserves a complex structure in $V$.
       Checking the Borel list of sphere transitive  linear groups, we conclude that
        $\vp(H) \simeq  H/ \ker(\vp)  \simeq  SU_m, U_m, Sp_{m/2}$ or $T^1 \cdot
        Sp_{m/2}$.\\
It remains  to check  that  the only possible case is  $\vp(H) =
U_m$ which is  clear for $m=1$.  For $m>1$, the   result  follows
from the following two lemmas.

\begin{lemma}\label{lemmaH/N}  
Let $H$ be    the  stability  subgroup  of  a  flag manifold $F
= G/H$. Then   there is no  normal  subgroup $N \subset H$  with
$H/N \simeq  SU_m,\, m>1$.
\end{lemma}

 \pf  If, by contradiction, a normal  subgroup $N \subset H$  such that
$H/N \simeq  SU_m$ exists, then there is an ideal  $\mathfrak{su}_m$ of $\mathfrak{h}$ corresponding to a  connected  component  of
    type $A_{m-1}$ in  the
    white   subdiagram  $\Pi_W$ of the painted Dynkin diagram
  $\Pi = \Pi_W \cup \Pi_B$ of the flag manifold $G/H$.

\noindent For the classical groups, we can always assume that  the  inclusion  of $A_{m-1}$  into $\Pi$  has  the  form

$$    \underset{\substack{\alpha_1} } {\circ} - \underset{\substack{\alpha_2}}{\circ} - \cdots -
 \underset{\substack{\alpha_{m-1}}}
{\circ}  - {\bullet} - \cdots
    $$

\smallskip

\noindent  Then the  ideal  $A_{m-1}= \mathfrak{su}_{m}$  is embedded into a subalgebra $A_{m} = \mathfrak{su}_{m+1}$
  and  the   subgroup  $H = N_G(\mathfrak{h}) \supset  N_{SU_{m+1}}(\mathfrak{su}_m)  =  U_m $ contains  $U_m$  as a normal  subgroup. Since $U_m = SU_m \times T^1/\mathbb{Z}_m$, we have $H/N \simeq SU_m/\mathbb{Z}_m$.
   This  implies  that
  there is  no normal  subgroup $N$  with  $H/N \simeq SU_m$, for $m > 1$, as required.\\

\smallskip

  Consider now  exceptional  groups. There   are only two   flag manifolds of  type  $F_4$, where  the  inclusion of the subalgebra $A_{m-1}= \mathfrak{su}_{m}$ is not as above:

$$
\underset{\substack{}}{\bullet} - \underset{\substack{}}{\bullet} \Rightarrow \underset{\substack{}}{\circ} - \underset{\substack{}}{\circ}  
$$

$$
\underset{\substack{}}{\circ} - \underset{\substack{}}{\circ} \Rightarrow \underset{\substack{}}{\bullet} - \underset{\substack{}}{\bullet}  \\
$$

\noindent  Indeed, in these cases the (white) subalgebra $\gh'=  A_2 = \mathfrak{su}_3$ is  embedded into $ \tilde{ \mathfrak{h}} = C_3 = \mathfrak{sp}_6$ for the first diagram and into
 $\tilde{ \mathfrak{h}} = B_3 = \mathfrak{so}_7$ for the second diagram. In both cases we have  $N_{\tilde{H}}(SU_3)  =  U_3 = SU_3 \times T^1/\mathbb{Z}_3$. Repeating  the above argument,   we  conclude that
   $H/N = SU_3/\mathbb{Z}_3$.

\smallskip

\noindent In the case when $G$ has  type $E_{\ell}, \, \ell = 6,7,8$, we can
always embed $A_{m-1}$ into $A_m$ or $D_m$  and  use   similar
arguments  with the exception of the cases

\bigskip

$$
\underset{}{\circ} - \underset{}{\circ} - \dotsb - \underset{}{\overset{\overset{\textstyle\bullet_{}}{\textstyle\vert}}{\circ}} \,-\, \underset{}{\circ} - \underset{}{\circ}  \\
$$

\noindent which correspond to the  flag manifold   $G/H = E_{\ell}/U_{\ell}, \,\, \ell=m = 6,7,8$.
 In this case  we have  $H = SU_{\ell} \cdot T^1$, that is, at the Lie algebra level, $\mathfrak{h} = \mathfrak{su}_l + i\mathbb{R}d$
  where  $\ad_d$ defines   a   gradation of depth 2 or 3
 of  the   complex Lie  algebra  $\mathfrak{e}_{\ell}^{\mathbb{C}}$. If we denote  by $V = \mathbb{C}^{\ell}$  the   standard $U_{\ell}$-module, then the  gradation  is  given  by (see \cite{G-O-V}, section 3.5) \\

 $$   \mathfrak{e}_{\ell}^{\mathbb{C}} = \mathfrak{u}_{\ell} + \mathfrak{g}_{-2} + \mathfrak{g}_{-1} + \mathfrak{g}_{1} + \mathfrak{g}_{2}$$
 \noindent for $\ell =6,7$, where $\mathfrak{g}_{-1} = \Lambda^3 V,\,  \mathfrak{g}_{-2} = \Lambda^6V,\, \mathfrak{g}_{1}= \Lambda^3 V^*, \, \mathfrak{g}_{-2}= \Lambda^4V^* $, and by
 $$   \mathfrak{e}_{8}^{\mathbb{C}} = \mathfrak{u}_{8} + \mathfrak{g}_{-3}+ \mathfrak{g}_{-2} + \mathfrak{g}_{-1} + \mathfrak{g}_{1} + \mathfrak{g}_{2} +\mathfrak{g}_{3} $$

\noindent where  $\mathfrak{g}_{\pm 1}, \mathfrak{g}_{\pm 2}$  are  defined
as above   and  $\mathfrak{g}_{-3}= V \otimes \Lambda^8 V,\,
 \mathfrak{g}_{3} =V^* \otimes  \Lambda^8 V^*$. The isotropy action of  the  semisimple  part  $SU_{\ell}$   of  the  stability  subgroup
 $H = SU_{\ell} \cdot T^1$  is the  standard action  on   forms.

\noindent Since the isotropy action  of $ i \cdot d$ on the    space  $\mathfrak{g}_k $ is  $i k \cdot \mathrm{id}$ and $\mathfrak{h} = \mathfrak{u}_{\ell} = \mathfrak{g}_0$, we have $ [d, \mathfrak{h}] = 0$, that is $ i \cdot d \in \mathfrak{z}(\mathfrak{h})$, so that $N= \exp (\mathbb{R}i\cdot d) = Z(H)$. By $Z(H) \cap SU_m = Z(U_m) \cap SU_m = \mathbb{Z}_{m}$ we have then

$$H/N = H/\exp (\mathbb{R}i\cdot d) = H/Z(H) = (H \cap SU_m)/(Z(H) \cap SU_m) = SU_m/\mathbb{Z}_{m}.$$

\noindent  The case of  the  group $G_2$ is  similar.\qed

\begin{lemma} 
There  is  no normal  subgroup $N \subset H$ with the  quotient
$H/N \simeq Sp_m$ or  $T^1\cdot Sp_m$ such that  the
 associated \C1 manifold $M = G\times_H V_{\vp}$
 has ordinary  regular orbits.
\end{lemma}

\pf  The   only   flag manifold  $G/H$ of a  classical Lie  group
which admits   a normal subgroup  $N \subset H$  with  indicated
quotient has the form $G/H = Sp_n/ (U_{n_1} \times \cdots \times U_{n_s}
\times T^k \times Sp_{m}) = Sp_n/ (N \times Sp_{m})$.
 The   associated \C1  manifold
$Sp_n \times_H V_{\vp} $  has    regular orbits given, for non-zero $v \in V_{\vp}$, by
$$Sp_n/L = Sp_n/N \times  (Sp_{m})_v  = Sp_n/ N \times Sp_{m-1} $$
\noindent which are  not ordinary  since
  $N_G(L) \setminus L \supset N_{Sp_m}(Sp_{m-1})\supset Sp_1$ and $Sp_1$ is not one-dimensional.

\smallskip

\noindent For the exceptional case, we have to consider the cases of $ F_4/ Sp_3 \cdot T^1 $  with  Dynkin diagram\\

$$
\underset{\substack{}}{\bullet} - \underset{\substack{}}{\circ} \Rightarrow \underset{\substack{}}{\circ} - \underset{\substack{}}{\circ}\\
$$

\noindent and $ F_4/Sp_2 \cdot T^2$  with Dynkin diagram\\

$$
\underset{\substack{}}{\bullet} - \underset{\substack{}}{\circ} \Rightarrow \underset{\substack{}}{\circ} - \underset{\substack{}}{\bullet} \\
$$

\noindent In both  cases, the associated  \C1 manifold   has non ordinary   regular orbits. For example, in the first case the corresponding  regular orbit  is
$F_4/(T^1 \cdot Sp_3)_v = F_4/Sp_2$, so that $N_G(L) \setminus L \supset N_{F_4}(Sp_2) \supset Sp_1$ and we conclude as in the classical case above. The case of $ F_4/Sp_2 \cdot T^2$ is  similar. This  finishes the proof
of the lemma. \qed

\subsubsection{Description of    admissible   homogeneous vector  bundles in terms of Dynkin diagrams}\label{SectionDescription}
The  following   proposition  describes the admissible homogeneous vector bundles
 $\pi: M_{\vp} = G \times_H V_{\vp}  \to   S_0 = G/H $, over a flag manifold
$(S_0 =G/H, J^S)$  in terms of painted  Dynkin diagrams   and      characters $\chi : T^k \to T^1$ of  the connected  center  $T^k$ of
 the   stabilizer   $H = H' \cdot T^k$.\\

\begin{proposition}\label{ pair}
  Let  $(S_0 =G/H = G/H' \cdot T^k, J^S)$ be  a  flag manifold   associated  with painted
    Dynkin  diagram  $\Pi = \Pi_B \cup \Pi_W$. Then an admissible homogeneous vector bundle $M_{\vp} = G\times_H V_{\vp}$
     is defined  by  a  pair   $(A_{m-1},  \chi)$  where
 $A_{m-1}$  is  a  connected component  of  $\Pi_W$ of  type $A_{m-1}$ (i.e. a string
of length $m-1$), 
and $\chi : T^k  \to T^1$ a character.
\end{proposition}

\pf    The  stability  subalgebra of $S_0$   admits  a direct sum  decomposition
  $$  \mathfrak{h}= \mathfrak{su}_m \oplus \mathfrak{n'} \oplus \mathrm{Z}(\mathfrak{h}) $$
  where  $\mathfrak{su}_m$ is   the ideal  associated  with   the   string $A_{m-1}$    and
 the   corresponding   decomposition  of    the   stability  subgroup is $H = SU_m \cdot N' \cdot T^k$. 
Then the extension  of  the   tautological representation of $SU_m$ in   a vector space  $V_{\vp } = \mathbb{C}^{m}$ by  a character
$\chi : T^k \to T^1 \cdot \mathrm{id}_{V_{\vp}} $   is  an  admissible    representation  $\vp : H \to  \vp(H) = U(V_{\vp}) \simeq  U_m$.
 The  converse   statement  is  also  clear.
\qed

 \subsection{Regular   submanifold    $M_{reg}$  of  $M_{\vp}$}
\subsubsection{ Invariant complex   structures  in the projective  space  $PV_{\vp} =H{[e_0]}=  U_m/U_{m-1}\times U(1)_0$}\label{Section621}
Let  $ \pi: M_{\vp} =  G \times_H V_{\vp} \to  S_0 =G/H$ be an admissible vector  bundle  over  a  flag manifold
$(S_0 =G/H, J^S)$  with      reductive  decomposition  $\mathfrak{g}= \mathfrak{h} + \mathfrak{p}$.\\

   We identify   $V_{\vp}$   with the   arithmetic complex vector   space   $\mathbb{C}^m$     with  the standard
  Hermitian    form  $<.,.>$    and    the   standard  basis
  $e_0 = (1, \cdots, 0)^T, \cdots , e_{m-1} = (0, \cdots, 0,1)^T$. Then  the  stability  subalgebra may be  written  as \\
  $\mathfrak{h} = \mathfrak{n} \oplus   \mathfrak{u}_m$    where   $\mathfrak{u}_m$ is  the   Lie   algebra of skew-Hermitian matrices.
   The  orbit  $\vp(H) e_0 =  H/H_{e_0}$ is     the  sphere with  the   reductive  decomposition of $\mathfrak{h} = \mathfrak{n}\oplus  \mathfrak{u}_m$   given  by
  $$ \mathfrak{h} = (\mathfrak{n} \oplus \mathfrak{u}_{m-1})   + ( \mathbb{R}I_0 + \mathfrak{q})$$
  where $  I_0 = \mathrm{diag}(i,0_{m-1} )$ and
   $$ \mathfrak{q}= \{  C_X :=
   \begin{pmatrix}
   0&-X^*\\
 X & 0
 \end{pmatrix} , \, X \in \mathbb{C}^{m-1} \} \simeq \mathbb{C}^{m-1}.
 $$
   The  diagonal Lie  algebra $\mathfrak{c}_{\mathfrak{u}_m}$    is  a Cartan   subalgebra of  $\mathfrak{u_m}$    and a basic  vector  $e_j$  is  a   weight  vector
   for   $\mathfrak{c}_{\mathfrak{u}_m}$     with     weight $\epsilon_j$  where
    $$\epsilon_j(\mathrm{diag}(x_0, x_1, \cdots, x_{m-1})) = x_j.$$
     The elementary matrices
   $E_{ij}  \in  \mathfrak{u}_m^{\mathbb{C}} = \mathfrak{gl}_m(\mathbb{C})$  are   root vectors    with  roots  $\alpha_{ij} = \epsilon_i - \epsilon_j$.
    The    reductive  decomposition     $\mathfrak{u}_m =  (\mathfrak{u}_{m-1} + \mathbb{R}I_0) + \mathfrak{q}$   of  the  projective  space  $\mathbb{C}P^{m-1} = U_m/U_{m-1}\times U(1)_0$
  defines a decomposition  of  the   root  system $R_{\mathfrak{u}_m}$ of $\mathfrak{u}_m^{\mathbb{C}}$ into  the  union  of white  roots $R^0_{\mathfrak{u}_m}  = \{ \epsilon_i - \epsilon_j, \, i,j >0\}$
   (  which  are   roots of   the  stability  subalgebra)  and  the complementary  black roots $R'_{\mathfrak{u}_m}  = \{ \pm \alpha_{0j}= \pm (\epsilon_0 - \epsilon_j),\, j>0 \}$.\\
   The multiplication  by  $\pm i$    in   $T_o \mathbb{C}P^{m-1} = \mathfrak{q} = \mathbb{C}^{m-1}$  defines  two (opposite)  invariant    complex  structures
    $\pm J_{\mathfrak{q}} $  which  define  two invariant  complex  structures   $\pm J^{\mathbb{C}P^{m-1}}$  on  $ \mathbb{C} P^{m-1}$.  They     correspond   to the   following painted Dynkin  diagrams:

    $$ \stackrel{\alpha_{01}}{\bullet}- \stackrel{\alpha_{12}}{\circ}- \cdots - \stackrel{\alpha_{m-1 m}}{\circ}\,\,\,\quad
 \stackrel{\alpha_{1 0}}{\bullet}- \stackrel{\alpha_{21}}{\circ}- \cdots - \stackrel{\alpha_{m m-1}}{\circ}  .$$
  Note    that $\alpha_{01}(\mathrm{id}_0) =1$  where   $\mathrm{id}_0 = -i I_0 = \mathrm{diag}(1,0, \cdots, 0) $. So   the $T$-Weyl  chamber   associated   with   the  complex  structure $\pm J^{\mathbb{C}P^{m-1}}$ is $\pm \mathbb{R}^+ \mathrm{id}_0$.\\

    Deleting  the   singular orbit  $S_0 = G/H$   which is    the   zero section  of $\pi$,  we get a   regular open   submanifold
    which  is    the  union  of   the    codimension one    regular  orbits  parametrized  by   $t >0$\\
    $ S_t = G (te_0)= G \times_H (te_0) = G/L $     where  $L = H_{e_0} =\mathrm{ker}(\vp) \cdot U_{m-1} $ is  the   stabilizer of   the point $e_0$.
    So  we  identify  $M_{reg}$    with
     $$   M_{reg} = \mathbb{R}^+ \times G/L.   $$
  The  projectivization  $G \times_H PV_{\vp}$  of  the  vector bundle  $\pi$   is a   homogeneous manifold \\
  $   F = G\{te_0\} = G \times_H [e_0] = G/K $   where   the    stabilizer $K = H_{[e_0]} = L \times U(1)_0$   with  $U(1)_0 = \mathrm{diag}(U(1), \mathrm{id_{m-1}}).$\\
We   will  assume  that the   regular orbit $S = G/L$ is  an ordinary manifold. Then $N_G(L) = K = L \times U(1)_0$  and   the  natural projection
$ G/L \to   F = G/K$  is    a principal $U(1)_0$-bundle over   the  flag manifold  $F$. The   restriction of  this projection  to  a   fibre  $V(x) = \pi^{-1}(x)$   is  the   standard  projection
  $S^{2m-1} = S_t \cap V(x) \to  PV(x)$  of  the sphere onto  the projective  space. To be  specific,   we  assume  that  $J_F = J_F^+$   corresponds  to  $\beta = \alpha_{12}$   such  that   the   system of   black  roots
$$\Pi^F_B = \{  \beta = \alpha_1, \beta_1, \cdots, \beta_{k}  \}$$
  and    the $T$-Weyl  chamber  is   defined   by  conditions
  $$C(J_F) = \{\beta >0, \beta_1>0, \cdots, \beta_k >0\}  \subset \mathfrak{t} = i Z(\mathfrak{k})  $$

 We have   the   following   standard  ($B$-orthogonal) decomposition
   $$  \mathfrak{g} = \mathfrak{l}+  \mathbb{R}Z^0_F + \mathfrak{m}= \mathfrak{l}+  \mathbb{R}Z^0_F + \mathfrak{q}  + \mathfrak{p} $$
where   $Z^0_F  \in  \mathbb{R} I_0 $ is      the    fundamental  vector, i.e.  the  vector of   $Z(\mathfrak{k}) = Z(\mathfrak{l}) \oplus  \mathbb{R}I_0$  orthogonal   to $\mathfrak{l}$ and normalized  by  the   condition   $B(Z^0_F, Z^0_F)=-1$.  Note  that   by  assumption,  its  centralizer  is  $\mathfrak{k}$.\\

\subsubsection{ Extension of  the   complex  structure  $J^S$    to   an invariant  complex  structure  in $F$}
Let  $\Pi = \Pi_W \cup \Pi_B $ be   the  PDD   associated   with  $(S_0,J^S)$   and  $\Pi_B = \{ \beta_1, \cdots, \beta_k\}$.
We  may  assume  that     the  $T$-Weyl  chamber $C(J^S)$ associated  with  the   complex  structure $J^S$ is   defined  by
$$C(J^S) =\{ \beta_1 >0, \cdots, \beta_k>0\} \subset  \mathfrak{t}_S = i Z(\mathfrak{h} ) = iZ(\mathfrak{n}) + \mathbb{R}\mathrm{id}_m .$$
 We  extend   the   complex  structure  $J_{\mathfrak{p}} = J^S|{\mathfrak{p}}$     to   a  $\mathrm{Ad}_K$-invariant  complex  structure  $J^{\pm}_{\mathfrak{m}}$
 on  $\mathfrak{m} = \mathfrak{q} + \mathfrak{p}$  choosing  one of  the   complex   structures    $\pm J_{\mathfrak{q}}$, described  above.  This   defines   two invariant   complex  structures $J^{\pm}_F$  in  the  flag manifold  $F = G/K$   which are consistent  with   the  complex  structure  $J^S$ in  $S_0$  such  that  the  natural projection  $F \to S_0$ is   holomorphic.
 The  PDD  of the   flag manifold   $(F, J^{\pm}_F)$ is  obtained   from  the  PDD of  $(S_0,J^S)$  by  painting in    black  one  of   the  end  roots  $ \beta = \alpha_1$ or  $\alpha_{m-1}$ of
  the  white  string  $A_{m-1} =\{ \alpha_1, \cdots, \alpha_{m-1}\}  \subset  \mathfrak{t}_F = i Z(\mathfrak{k}) = \mathfrak{t}_S + \mathbb{R}Z^0,\,\,  Z_0 = -i Z^0_F $.\\
Since   $\beta(Z^0) \neq 0$,  changing  $Z^0$ to $- Z^0$ if necessary,  we may  assume  that  $\beta(Z^0) >0$.\\

\subsubsection{ Decomposition    $TM_{\vp} =  T^hM_{\vp} + T^v M_{\vp}$ of  the tangent  bundle }

 The  decomposition  $\mathfrak{g}= \mathfrak{h}+ \mathfrak{p}$
defines  a  $G$-invariant  principal  connection  in  the principal bundle  $G \to  S_0 = G/H$  such  that   the   decomposition
 of   the tangent   bundle  $TG$ into  horizontal    subbundle  $T^hG$ and  vertical  subbundle  is  given  by
   $$ T_aG   = T^h_a G + T^v_aG =  a \mathfrak{p}  + a \mathfrak{h}= (L_a)_*\mathfrak{p}+ (L_a)_*\mathfrak{h} .$$
It  defines a similar   decomposition  of  the  tangent  bundle  $T M_{\vp}$  given  by
$$   T_{[a,v]}M_{\vp} = T^h_{[a,v]} M_{\vp}  +  T^v_{[a,v]} M_{\vp}=   T^h_aG + T_vV_{\vp} = a \mathfrak{p} + V_{\vp} . $$
 In particular,   along  the  radial  line  $\mathbb{R} e_0 \in  \pi^{-1}(eH)\subset  M_{\vp}$    the  tangent  space  $T_{te_0}M_{reg}$ can be  written   as
 $$ T_{te_0}M_{\vp} = T^h_{te_0} M_{reg} +  T^v_{t e_0}M_{\vp} =  \mathfrak{p} + V_{\vp}.  $$

 \begin{proposition}\label{propVM}    
Let $g_V$  be   a  $\vp(H) = U_m$-invariant  metric   in  $V_{\vp}$   and  $g^t_{\mathfrak{p}},\,  t \in \mathbb{R}$ a 1-parameter  family  of   $\mathrm{Ad}_K$-invariant  metrics in
  $\mathfrak{p}$.\\
 Then  the  metric  $g_{t e_0} = g^t_{\mathfrak{p}}  \oplus g_V$ in  $\mathfrak{p} + V_{\vp} = T_{te_0} M_{reg} $
 is  extended   by left translations  from  $G$  to a   smooth invariant metric in  $M_{\vp}$.
 \end{proposition}

\pf Indeed, the metric  $g_{\mathfrak{p}}$ defines   an invariant  metric  in  $S_0 = G/H$  which induces   the $G$-invariant metric $g^h$  in  the   horizontal  distribution  $T^h M_{\vp}$  via  the isomorphism  $\pi_* : T^h_x M_{\vp} \to  T_{\pi(x)} S_0$.  The  metric  $g_V$ in  the   fibre $V_{\vp} = \pi^{-1}(eH)$ induces   an  invariant metric $g^v$ in  the   vertical distribution $T^vM_{\vp}$.
Then   $g = g^h \oplus g^v$ is   an invariant Riemannian metric in  $M_{\vp}$.
\qed

\subsubsection{Description  of  the   character}

 For simplicity,  we will  assume that  $G$ acts effectively on
$S_0$. Then $G$ has no center  and we may  identify
  $G$   with  the  adjoint  group $\mathrm{Ad}_G$
 and $H$  with  the   adjoint  subgroup $\mathrm{Ad}_H \subset
 \mathrm{Ad}_G$.
  Then the group of characters $\mathcal{X}(T^k)$ is identified
  with the  lattice  $Q_T = \mathrm{span}_{\mathbb{Z}}R_T \subset \gt_H^*$  of $T$-roots as
  follows.\\
We denote  by
$$   Q_T^* := \{ h \in \gt,\,  \bar \beta(h) \in \bZ , \bar \beta \in R_T \}=
\mathrm{span}(h_1, \cdots, h_k) \subset \gt_H := i Z(\mathfrak{h})
$$
the  lattice in $\gt_H =i Z(\mathfrak{h})$  dual  to  $Q_T$  with
basis  $(h_i,\, i=1, \cdots, k)$ dual
 to the  basis  $\{ \bar \beta_j = \kappa(\beta_j),\,  \beta_j \in \Pi_B \} $ of  $Q_T$,
  $\bar \beta_j(h_k)= \delta_{jk}$.

 \noindent Then $2 \pi Q_T^*$ is  the kernel of   the  exponential map
  $$   \exp \circ i : \gt \to T^k,\,  h \mapsto \exp\,ih $$
  and  we  may  identify $T^k$  with   $\mathfrak{t}_H/ 2\pi ( Q_T)^*$.

\noindent A weight
$$\Lambda = \sum_{j=1}^k p_j \bar \beta_j \in  Q_T =
\mathcal{X}(T^k):  h = \sum_{j=1}^k x_j h_j \mapsto  \Lambda(h) = \sum_{j=1}^k  p_jx_j $$
   defines a  character
  $$ \chi = \chi_{\Lambda} : T^k \to \bC^*,\,\,  \exp ih \mapsto e^{2\pi i \Lambda(h)} = e^{2\pi  i p_j x_j}
.$$

\noindent (The reader is referred for example to \cite{Sn} for the correspondence between homogeneous line bundles and characters).

\medskip

\subsection{Invariant K\"ahler  structures  on  $M_{\vp}$}

  Let  $(S_0 =G/H, J^S)$ be  a  flag manifold  with  reductive  decomposition   $\mathfrak{g}= \mathfrak{h}+ \mathfrak{p}$  associated  with painted
    Dynkin  diagram  $\Pi = \Pi_W \cup \Pi_B,\, \Pi_B = \{  \beta_1, \cdots, \beta_k\}$   and   the   complex  structure $J^S$   associated  with  the  $T$-Weyl  chamber
    $C(J^S) = \{  \beta_1 >0, \cdots, \beta_k >0\}  \subset  \mathfrak{t}_S = i Z(\mathfrak{h})$.
    Let $M_{\vp} = G\times_H V_{\vp}, \,  V_{\vp} = \mathbb{C}^{m-1}$  be   the admissible homogeneous vector bundle   associated   with  a  pair   $(A_{m-1},  \chi)$  where
 $A_{m-1}$  is  a  connected component  of  $\Pi_W$ of  type $A_{m-1}$ (i.e. a white  string
of length $m-1$     and $\chi : T^{\ell} =Z^{con}(H) \to   T^1 $ a character.\\
We  assume  that   the regular orbit  $S_t = G(te_0) = G/L$ is ordinary. Then   the  projectivization $F = G\times_H PV_{\vp}  = G/K = G/L\times U(1)$  is  a   flag manifold.
 Let  $J^F$  be    an  extension of  the  complex   structure $J^S$  to  an invariant  complex  structure of $F$    defined by  extension  of    the    complex  structure  $J_{\mathfrak{p}}=
 J^S|_{\mathfrak{p}}$  to   a  complex  structure  $J_{\mathfrak{m}} = J_{\mathfrak{q}} \oplus  J_{\mathfrak{p}}$ on  $\mathfrak{m = \mathfrak{q}+ \mathfrak{p}} = T_{eK}F$  such  that
 the  PDD  of  $(F= G/K,J^F)$  is  obtained  by  painting in  black   the   end   root  $\beta$ of  the  string  $A_{m-1}$.\\
  The standard   decomposition  associated  with the C1 manifold  $M_{reg}$  may be  written as
  $$\mathfrak{g}= \mathfrak{l}+ \mathbb{R}Z^0_F + \mathfrak{m} = (\mathfrak{n}+ \mathfrak{su}_{m-1}+ \mathbb{R}I_{m-1}) + (\mathfrak{q} + \mathfrak{p})$$
  where $\mathfrak{n} = \mathrm{ker}(\vp)$,  $\mathfrak{k} = \mathfrak{l} + \mathbb{R}Z^0_F$  and  $Z^0_F$ is  the  fundamental vector   which  is identified  with  the   fundamental vector  field on $M_{reg}$ whose  restriction  to    a regular orbit  $S_t$ is  the fundamental   vector  field
 of  the principal $U(1)$ -bundle  $S_t =  G/L \to  F = G/K = G/L \times U(1). $\\

 The main  result of  this  section is  the  following  theorem.

\begin{theorem}\label{KaehlerStructonM{vp}}
Let $M_{\vp} = G\times_H V_{\vp}$ be as above a standard C1 manifold having as singular orbit the flag manifold $(S_0 =G/H, J^S)$, where the complex  structure $J^S$ is associated  with  the  $T$-Weyl  chamber
    $C(J^S) = \{\beta_1>0, \cdots , \beta_k >0\}  \subset i Z(\mathfrak{h})$ and let $J^F$  be an  extension of  the  complex   structure $J^S$  to  an invariant  complex  structure of $F = G\times_H PV_{\vp}  = G/K$ associated with the $T$-Weyl chamber $C(J^F) = \{ \beta>0, \beta_1>0, \cdots , \beta_k >0\}$.

\noindent Let $Z^0_F = \kappa I_0$ be  the  fundamental vector, where $\kappa^{-1} = (-B(I_0, I_0))^{\frac{1}{2}}$ and $I_0$ is defined in Section \ref{Section621}, and let $\theta^0$ be the    dual to  $Z^0_F$ invariant   1-form on the regular part $M_{reg}$ of $M_{\vp}$.
 We may  assume  that   $\beta(Z^0) >0$  where $Z^0 = -i Z^0_F \in i Z(\mathfrak{k})$.

\noindent Let  $Z_0 \in  C(J^S)$  be    a  vector   from  the  $T$-Weyl  chamber  $C(J^S)$ which is   the   face of   the Weyl  chamber   $C(J^F) = \{ \beta =0\} $ .
   Then  a segment  in $C(J^F) $   with  a parametrization  $ Z_0 +  f(t) Z^0, \, { \dot f}(t) >0$   defines  a K\"ahler  metric in  $M_{reg}$ given  by
   (see Proposition \ref{KaehlerstructonM{reg}})

$$ g_{reg} = dt^2 + (\dot f \theta^0)^2 + \pi_F^* g_0 + f(t) \pi_F^*g^0,
$$
where  $\pi_F : M_{reg} \to  F$ is   the  natural projection,  $g_0 =  \omega_{Z_0}\circ J^F $ is  a  symmetric  bilinear  form  on $F$ ( which is  the pull back of  an  the invariant K\"ahler  metric  on $S_0$ ) and  $g^0 = -\omega_{Z^0}\circ J^F$ is    a   symmetric
bilinear   form on $F$.\\
The K\"ahler  structure   smoothly   extends  to   a  geodesically   complete  invariant  K\"ahler  structure  on  $M_{\vp}$ if  and only if   the  function  $f(t)$ is  extended  to a  smooth even   function on  $\mathbb{R}$ such that $Z_0 + f(t) Z^0 \in C(J^F)$ and satisfies  the   following  Verdiani  conditions :
$$f(0) = \dot f(0)=0,\,  \ddot f(0) =\kappa$$
\end{theorem}

\pf In the notations used before the statement, let $X \in  \mathfrak{q}, \, Y \in \mathfrak{m}$, then
 $$g_o (X,Y) =  d (B \circ Z_0)(X,J_{\mathfrak{m}}Y)= - B(Z_0, [X, J_{\mathfrak{m}}Y]) = $$
$$= - B([Z_0, X], J_{\mathfrak{m}}Y])=0.$$
This  shows   that  $\pi_S^*g_0$ is    a metric  in   the    horizontal    subbundle   $T^hM_{\vp}$.
 Similarly, for $X \in  \mathfrak{q}, \, Y \in \mathfrak{p}$ ,    we get
 $$ g^0 (X,Y) = - B(Z^0_F, [X, J_{\mathfrak{p}}Y] ) = -B([Z^0_F, X], J_{\mathfrak{p}}Y)=0  $$
since  $[Z^0_F, X] \in  \mathfrak{q}$ and  $J_{\mathfrak{p}}Y \in  \mathfrak{p}$.
This   shows  that  $g_{reg}(T^h M_{reg}, T^v M_{reg})  =0$.
 The    horizontal part $g_{reg}^h$  of  the  metric $g_{reg}$ at   a point  $t e_0$ can be   written  as
 $$  g_{reg}^h = g_{reg}|_{\mathfrak{p}} = \pi^*_F g_0 + f \pi^*_F g^0 .$$
Under the assumptions of the  theorem, it is  extended  to  a  smooth metric   in $T^hM_{reg}$.
The    vertical part  is
$g_{reg}^v = g_{reg}|_{V}   = dt^2 + (\dot f \theta^0)^2 + f \pi^*_F g^0|_{\mathfrak{q}}. $
  Since  $Z^0_F$   and   $\mathfrak{q}$ belong  to  $\mathfrak{u}_m \subset \mathfrak{h}$, using  calculation in
  $\mathfrak{u}_m$  we   get for  $X,Y \in \mathfrak{q}$
   $g^v_{reg}(X,Y) =  -d(B \circ Z^0_F)(X,iY) = B(Z^0_F, [X,iY])=  B(Z^0_F,  X^*iY+ iY^*X)=-2 <X,Y>B(Z^0_F, I_0)= \frac{2}{\kappa}<X,Y> $
  where   $<X,Y> =  Im (X^*Y) = \frac12(X^*Y + Y^*X)$ is  the    standard metric in  $\mathfrak{q} = \mathbb{C}^{m-1}$.

	\noindent Now, let $g_{eucl} = dt^2 + t^2 \eta^2 + g_{eucl}|_{\mathfrak{q}}$ be the flat euclidean metric on $V$, where $\eta$ denotes the 1-form dual to the vector field generated by $I_0$. Since $\theta^0 = \frac{1}{\kappa} \eta$ we can rewrite $g_{reg}|_{V}   = dt^2 + (\frac{\dot f}{\kappa t})^2 t^2 \eta^2 + f \pi^*_F g^0|_{\mathfrak{q}}$ and apply Verdiani   criterion (Theorem 1 in \cite{V}) which says that $g^v_{reg}|_{V_{\vp} \setminus \{0\}}$ is  extended  to $V_{\vp}$  if  and  only if $f(t)$ satisfies the stated conditions, i.e. is a smooth even function on $\mathbb{R}$ with $\dot{f}(0)=0$, $\ddot{f}(0)= \kappa$.  Then  by Proposition \ref{propVM} the metric $g_{reg}$ is  extended  to  a
      Riemannian metric  $g_M$ on  $M_{\varphi}$. Since    this   metric  is  K\"ahler  on   $M_{reg}$,  the   corresponding  complex  structure
      $J_{reg}$ is parallel. It is  clear   that  $J_{reg}$  is  extended  to  a parallel  complex  structure  on  $M_{\vp}$. Hence  the metric
      $g_M$ is  K\"ahler.\\

   \noindent The claim    about  the completeness of the metric on  $M_{\vp}$ follows from the following

\begin{lemma}
A metric of the form $dt^2 + g_t$ on a manifold $M = \bR \times N$, where $g_t$, $t \in \bR$ is a family of metrics on the compact manifold $N$, is complete.
\end{lemma}

\noindent {\bf Proof of the Lemma: } By Hopf-Rinow theorem, a metric is complete if and only if the closed balls are compact. This is true under our assumptions since any closed ball of radius $r$ in $M$ is contained in the compact set $[-r, r] \times N$. \qed

\begin{definition}
An interval, together with a parametrization as described in Theorem \ref{KaehlerStructonM{vp}}, is called admissible.
\end{definition}

\begin{remark}
 Notice that  for any $d > 0$, an admissible parametrization $f: (0, + \infty) \rightarrow (0, d)$  exists, take for example the function

$$f(t) = d (1 - e^{-\frac{\kappa}{2 d} t^2}).$$

\noindent This shows that, without additional conditions, a bounded segment in the Weyl chamber can define a complete K\"ahler metric (compare with the case of K\"ahler-Einstein metrics in the next section).
\end{remark}

\section{K\"ahler-Einstein metrics on  standard \C1 manifolds }

In this section we give necessary and sufficient conditions for the existence of (complete and non-complete) K\"ahler-Einstein metrics on standard C1 manifolds.
\smallskip

\noindent Our first result is the following

\begin{theorem}\label{maintheoronKE}
Let $M$ be a standard \C1 manifold, i.e. (see Theorem \ref{C1  complex manifolds}) the total space of an admissible bundle $M_{\vp} = G \times_H V_{\vp} \rightarrow S_0$ over the singular orbit $(S_0 = G/H, J^S)$.

\noindent Let $(F = G \times_H P V_{\vp} = G/K, J^F)$ be the flag manifold associated with regular orbits, $Z^{Kos} \in C(J^F)$ be the Koszul vector which defines  the  invariant K\"ahler-Einstein metric on $(F,J^F)$ associated  with the invariant complex  structure $J^F$ and $(Z_0 Z_d) \subseteq C = C(J^F)$ an interval in  the  $T$-Weyl  chamber  $C(J^F)$ which represents a standard invariant K\"ahler structure on $M_{\vp}$ (Theorem \ref{KaehlerStructonM{vp}}).

\noindent Then,    the  interval $(Z_0 Z_d)$    with  a parametrization  $f(t)$    defines   a K\"ahler-Einstein structure with Einstein constant $\lambda$ if  and only if
 the   vectors    $Z^{Kos},  Z_0, Z^0$      are  related  by

\be\label{algCond}
Z^{Kos} = \lambda Z_0 + \kappa m Z^0
\ee

 \noindent where $m = \dim(V_{\vp})$, $k$ is defined by $Z^0_F = \kappa I_o$ (Theorem \ref{KaehlerStructonM{vp}})
  and  $f(t)$ is      the   solution of  the   equation

   \be \label{EinstEq}
   \ddot f(t) + \frac12 A(f) \dot f^2 + \lambda f =\kappa m
   \ee

\noindent with  the initial  conditions\\
   $\lim_{t \to 0}f(t) = \lim_{t \to 0}\dot f(t)=0,\, \lim_{t \to 0}\ddot f(t)= \kappa$

\noindent where $A(f) = \sum_{\alpha \in R_{\mathfrak{m}}^+} \frac{\alpha(Z^0)}{\alpha(Z_0) +f \alpha(Z^0)}$, being $R_{\mathfrak{m}}^+$ the set of the positive black roots of $G/K$ (see formulas (\ref{Einstein1}) and (\ref{Einstein2})).

\noindent Moreover, the K\"ahler-Einstein metric can be extended to a complete metric if and only if $\lambda \leq 0$, and in this case the segment extends to a ray $Z_0 + \bR^+ Z^0$ in $C(J^F)$.
\end{theorem}

\pf

\noindent The calculations made in Section  \ref{einstequation2} show that the K\"ahler metric determined by the interval $(Z_0 Z_d)$ together with a parametrization $Z(t) = Z_0 + f(t)Z^0$ satisfies the Einstein condition if and only if

\be
Z^{Kos} = \lambda Z_0 + c Z^0
\ee

 \be \label{EinstEq}
\ddot f(t) + \frac12 A(f) \dot f^2 + \lambda f =c
 \ee

\noindent for some constant $c$. So, we must show that $c = \kappa m$.

\noindent In order to do that, recall that in Theorem \ref{KaehlerStructonM{vp}} it was shown that the metric extends to the singular orbit if and only if the function $f$ satisfies the initial conditions
 $\lim_{t \to 0}f(t) = \lim_{t \to 0}\dot f(t)=0,\, \lim_{t \to 0}\ddot f(t)= \kappa$. This implies that $f(t) = \frac{\kappa}{2}t^2 + O(t^3)$, and then
$$\frac12 A(f) \dot f^2 = \sum_{\alpha \in R_{\mathfrak{m}}^+} \frac{\alpha(Z^0)(\kappa^2 t^2 + O(t^3))}{\alpha(Z_0) + \alpha(Z^0)(\frac{\kappa}{2}t^2 + O(t^3))}$$
\noindent when $t \rightarrow 0$ tends to $0$ if $\alpha(Z_0) \neq 0$ and to $2 \kappa$ if $\alpha(Z_0) = 0$.
\noindent Since $Z_0 \in C_{\beta} = \{ \beta = 0, \beta_1 >0, \dots, \beta_k > 0 \}$ where $\{ \beta_1, \dots, \beta_k \}$ (resp. $\{ \beta, \beta_1, \dots, \beta_k \}$) is the set of black roots in the Dynkin diagram of $G/H$ (resp. $G/K$), then a positive black root $\alpha \in R_{\mathfrak{m}}^+$ is a black root of $G/H$ if and only if $\alpha(Z_0) \neq 0$ (recall that, since $Z_0 \in Z(\mathfrak{k})$ then every white root vanishes on $Z_0$).

\noindent In other words, the number of roots $\alpha \in R_{\mathfrak{m}}^+$ for which $\alpha(Z_0) = 0$ equals the number of positive black roots of $G/K$ minus the number of positive black roots of $G/H$, i.e. equals $dim_{\bC}(G/K) - dim_{\bC}(G/H)$ which, by $G/K =G \times_H P V_{\vp}$, is equal to $m-1$.
\noindent It follows then that $\frac12 A(f) \dot f^2 \rightarrow \kappa (m-1)$ for $t \rightarrow 0$, which, combined with the other initial conditions, implies that
$$\ddot f(t) + \frac12 A(f) \dot f^2 + \lambda f \rightarrow \kappa + \kappa (m-1) = \kappa m.$$
\noindent This  shows that $c = \kappa m$.

\noindent Notice that $f(t)$ extends to a smooth even function. In fact it follows by a straight calculation using equation (\ref{EinstEq}) that, under the given initial conditions, $\lim_{t \to 0} f^{(3)}(t)=0$, which shows that $f$ extends to a $C^3$ function invariant by reflection at $0$. Then, it gives rise to a $C^2$-Einstein metric and we can apply a result by DeTurck and Kazdan (see, for example, \cite{B}) to conclude that $f$ is $C^{\infty}$.

\medskip

\noindent In order to end the proof of the theorem, we need to prove the following

\begin{lemma}\label{theoremDescribesKE}
If the condition (\ref{algCond}) is fulfilled, then the function $f(t)$ parametrizing the segment $(Z_0 Z_d)$ which gives the K\"ahler-Einstein metric is the inverse to the function

\be\label{t(f)}
t(f) = \int_{0}^{f} \sqrt{\frac{P(s)}{2 \int_0^s (c - \lambda v) P(v) dv}} ds
\ee

\noindent where $P$ is the polynomial defined by $P(x) = \Pi_{\alpha \in R_{\mathfrak{m}}^+} (\alpha(Z_0) + x \ \alpha(Z^0))$.
\end{lemma}

\pf
The proof is based on the fact that, if $f$ satisfies the ordinary differential equation $\ddot f(t) + \frac12 A(f) \dot f^2 + \lambda f =c$, where $A(f)$ is any function of $f$, then

\be\label{LemmaODE}
(f')^2 = e^{-\int_0^f A(v) dv} \left( \int_0^f 2(c - \lambda v) e^{\int_0^v A(s) ds} dv + D \right)
\ee

\noindent Indeed, by the substitution $p(f) = f'$ we get $f''(t) = p'(f) f' = p' p$, so that the equation can be rewritten
$$p' p + \frac12 A(f) p^2 + \lambda f =c$$
\noindent that is, by setting $u(f) = p^2(f)$,
$$u' = -A(f) u + 2 (c - \lambda f).$$

\noindent Then (\ref{LemmaODE}) follows by using the general formula to solve a linear first order ODE $x' = Px+ Q$, that is $x(t) = e^{\int P(t) dt} \left[ \int Q(t)  e^{-\int P(t)dt} dt + D \right]$.

\noindent In particular, in our case $A(f) = \sum_{\alpha \in R_{\mathfrak{m}}^+} \frac{\alpha(Z^0)}{\alpha(Z_0) +f \alpha(Z^0)}$, we have $e^{\int A(f) df} = \Pi_{\alpha \in R_{\mathfrak{m}}^+} (\alpha(Z_0) +f \alpha(Z^0)) = P(f)$ which, replaced into (\ref{LemmaODE}) yields

\be\label{LemmaODE2}
(f')^2 = \frac{1}{P(f)} \left( \int_0^f 2(c - \lambda v) P(v) dv + D \right)
\ee

\noindent from which formula (\ref{t(f)}) follows by extracting the square root, inverting and integrating, and by taking into account that for $t \rightarrow 0$ we have $f(t) \rightarrow 0$ and $f'(t) \rightarrow 0$.
\qed

\bigskip

\noindent Now, using formula (\ref{t(f)}) we are ready to end the proof of Theorem \ref{maintheoronKE}.

\smallskip

\noindent If $\lambda > 0$, then the metric is not complete by Myers' theorem, which implies that a complete Kahler-Einstein manifold with positive Einstein constant is compact (and then in the C1 case should have two singular orbits).

\noindent Let then $\lambda \leq 0$. Then, by (\ref{algCond}) we have $\kappa m Z^0 = Z^{Kos} - \lambda Z_0$ and since $Z^{Kos}, Z_0$ belong to $C(J^F)$ then also $Z^0 \in C(J^F)$.

\noindent This implies that the polynomial $P(x) = \Pi_{\alpha \in R_{\mathfrak{m}}^+} (\alpha(Z_0) + x \ \alpha(Z^0))$ is strictly positive for $x > 0$: indeed,  by definition, a root is black if and only if in its decomposition as sum of the basis roots there is one of those corresponding to the black vertices of the painted  Dynkin diagram, that is $\beta, \beta_1, \dots, \beta_p$. Since for any white root $\gamma$ one has $\gamma(Z_0) = \gamma(Z^0) = 0$ (since $Z_0, Z^0 \in Z(\mathfrak{k})$), the factors $\alpha(Z_0) + x \ \alpha(Z^0)$ of the polynomial $P(x)$ are positive for $x >0$ if and only if $\beta(Z_0) + x \ \beta(Z^0) = \beta(Z_0 + x Z^0)$, $\beta_j(Z_0) + x \ \beta_j(Z^0) = \beta_j(Z_0 + x Z^0)$, $j=1, \dots, p$, are positive, which is true by definition of $C(J^F)$.

\noindent From this and from $\lambda \leq 0$, it follows that the polynomial $Q(s) = 2 \int_0^s (\kappa m - \lambda v) P(v) dv$ appearing in the formula (\ref{t(f)}) has first derivative $Q'(s) = 2 (\kappa m - \lambda s) P(s)$ always positive for $s > 0$, so it is strictly increasing and (being $Q(0)=0$) we have $Q(s) > 0$ for any $s$.

\noindent Then the integrand $\sqrt{\frac{P(s)}{2 \int_0^s (\kappa m - \lambda v) P(v) dv}}$ in (\ref{t(f)}) is always well-definite and obviously positive. Moreover, being the square root of the ratio between one polynomial of degree $N$ and a polynomial of degree $N+2$ (both with no positive real roots), it goes to infinity like $1/s$, and then if $f \rightarrow \infty$ then one has $t \rightarrow \infty$: this shows that $f$ does not blow in finite time, and then it is defined on $(0, + \infty)$, which, together with the fact that $Z_0 + f(t) Z^0 \in C(J^F)$ (we have already observed above that $Z^0 \in C(J^F)$) proves that the metric is complete.

\noindent Now, in order to prove that $Z_0 + fZ^0$ is actually a ray in $C(J^F)$, we observe that the function $f$ does not stay bounded.

\noindent Indeed, if so, then either there would be a point $t_0 > 0$ for which if $f'(t_0) = 0$, or $f'(t) \rightarrow 0$ for $t \rightarrow \infty$. In the first case, let $t_0$ be the first positive value for which if $f'(t_0) = 0$: by the initial conditions, it should be $f''(t_0) \leq 0$. But, from equation (\ref{EinstEq}), one gets $f''(t_0) = c - \lambda f(t_0)$ and then, since $c = \kappa m > 0$ and $\lambda \leq 0$, we have $f''(t_0) > 0$, a contradiction. In the second case, we can conclude by the same argument since, from equation (\ref{EinstEq}),

$$f'' = c - \lambda f - \frac{1}{2} A(f) f'^2 \rightarrow c - \lambda f.$$

\qed

\bigskip

\noindent Using Theorem \ref{maintheoronKE} and Lemma \ref{theoremDescribesKE}, together with the description of standard C1 manifolds in terms of Dynkin diagrams and characters given in Section 6, it is possible to find explicit conditions (at least when $G$ is a classical group) for the existence and completeness of K\"ahler-Einstein standard C1 manifolds having a given flag manifold $S_0 = G/H$ as (the only) singular orbit, for any sign of the Einstein constant. More precisely, one has

\begin{proposition}\label{prop triple2}
  Let $G$ a simply connected Lie group with Lie algebra $\mathfrak{g}$ equal to one of the classical Lie algebras $su_n, sp_{2n}, so_{2n}, so_{2n+1}$, and let $(S_0 =G/H, J^S)$ be  a  flag manifold   associated  with painted
    Dynkin  diagram  $\Pi = \Pi^H_B \cup \Pi^H_W$ (possibly consisting of more connected components) which begins with a white $A_{m-1}$ string.
		
\noindent Let $G/K$ be the flag manifold obtained by painting in black one of the end roots of the $A_{m-1}$ string, and let $n_1, \dots, n_p$ be the coefficients of the fundamental weights associated to the roots in $\Pi^H_B$ in the decomposition of the Koszul form of $G/K$. 	
		
		Then,

\begin{enumerate}
\item[(i)] there exists a K\"ahler-Einstein standard C1 manifold $M$ (having $S_0$ as only singular orbit) with Einstein constant $\lambda=0$ if and only if $n_1, n_2, \dots, n_p$ (resp. $n_1+1, n_2, \dots, n_p$) are divisible by $m$ if the $A_{m-1}$ string is a connected component of the diagram (resp. otherwise). In particular, if $m=1$ this condition is trivially fulfilled.

\item[(ii)] for any $\lambda \neq 0$, there always exists a K\"ahler-Einstein standard C1 manifold $M$ (having $S_0$ as only singular orbit) with Einstein constant $\lambda$.

\end{enumerate}
\end{proposition}

\noindent If the above conditions are satisfied then, in accordance with above Theorem \ref{maintheoronKE}, the Kahler-Einstein metric can be chosen to be complete for $\lambda \leq 0$, while for $\lambda > 0$ the metric is never complete. Notice that these results are obviously consistent with Myers' theorem and with \cite{D-W}.

\noindent We give a proof of Proposition \ref{prop triple2} and illustrate it with some examples in the second part of this paper.

\end{document}